\documentclass[a4paper,11pt]{amsart}

\usepackage{amsmath, amssymb, amsthm}
\usepackage{verbatim}
\usepackage{graphicx}
\usepackage[top=3cm, bottom=3cm, left=2.5cm, right=2.5cm]{geometry}
\usepackage[colorlinks=true, linkcolor=blue, citecolor = blue]{hyperref}
\usepackage{mathrsfs}

\usepackage{mathrsfs}
\usepackage{psfrag}
\usepackage{colortbl}
\usepackage[normalem]{ulem}
\usepackage[percent]{overpic}
\usepackage{diagbox}

\newtheorem{remark}{Remark}[section]
\newtheorem{definition}{Definition}[section]
\newtheorem{theorem}{Theorem}[section]
\newtheorem{lemma}{Lemma}[section]

\newcommand{\R}{\mathbb{R}}
\newcommand{\eps}{\varepsilon}
\newcommand{\deff}{\buildrel {\rm def}\over =}

\newcommand{\jsn}{{\tt{\,sn}}}
\newcommand{\jcn}{{\tt{\,cn}}}
\newcommand{\jdn}{{\tt{\,dn}}}
\newcommand{\omg}{\omega}
\newcommand{\Omg}{\Omega}
\newcommand{\etat}{{\eta}}
\newcommand{\deltat}{{\delta}}
\newcommand{\sigmat}{{\sigma}}
\newcommand{\deltah}{{\delta}}

\newcommand{\deltahat}{\hat{\delta}}
\newcommand{\lambdahat}{\hat{\lambda}}
\newcommand{\gammahat}{\hat{\gamma}}
\newcommand{\expt}{\mathbb{E}}

\newcommand{\frpt}[2]{\frac{\partial #1}{\partial #2}}
\newcommand{\frptsec}[2]{\frac{\partial^2 #1}{\partial^2 #2}}

\newcommand{\psih}{\widehat{\psi}}
\newcommand{\la}{\langle}
\newcommand{\ra}{\rangle}
\newcommand{\mfF}{\mathfrak{F}}
\newcommand{\mfG}{\mathfrak{G}}
\newcommand{\mfB}{\mathfrak{B}}
\newcommand{\Aavg}{\mathbb{A}}
\newcommand{\bbK}{\mathbb{K}}
\newcommand{\mfT}{\mathfrak{T}}
\newcommand{\mfc}{\mathfrak{c}}
\newcommand{\mfs}{\mathfrak{s}}
\newcommand{\msV}{\mathscr{V}}
\newcommand{\mfh}{\mathfrak{h}}
\newcommand{\stopt}{\mathfrak{e}}
\newcommand{\stopteps}{\stopt^\eps}
\newcommand{\msB}{\mathscr{B}}

\newcommand{\hbarst}{\mfh^0}

\title[Escape from a resonance zone]{Random perturbations of a periodically driven nonlinear oscillator: Escape from a resonance zone}
\author{Nishanth Lingala}
\author{N. Sri Namachchivaya}
\author{Ilya Pavlyukevich}



\begin{document}

%#####################################################################
%#####################################################################
%#####################################################################
%#####################################################################

\begin{abstract} 
The phase space for the periodically driven nonlinear oscillator
consists of many resonance zones. Let the strength of periodic excitation and the strength of the damping be indexed by a small parameter $\eps$. It is well known that, as $\eps\to 0$, the measure of the set of initial conditions which lead to capture in a resonance zone goes to zero. In this paper we study the effect of weak noise on the escape from a resonance zone.
\end{abstract}

\maketitle

%##################
%##################
%##################
\section{Introduction}
\label{sec:intro}
\noindent 
The recent surge of research articles in energy harvesting focuses on the ``cantilever beam" devices which are used to convert small amplitude mechanical vibration into electrical energy that could be used for electronic devices with low power requirements (see \cite{DaqMasSurv}). Prototypical beam type nonlinear energy harvesting models contain double well potentials, external or parametric periodic forcing terms, damping and ambient broadband additive noise terms. For example, \cite{McInnes2008,Cartmell2014} considers devices which can be modeled by
\begin{equation}\label{E:Cartmell}
\ddot q_t + \delta \dot q_t  - \mu (1 + \eta\cos(\nu t) ) q_t +\gamma q_t^3=  \sigma \,\xi(t) +\alpha \cos(\nu t) ,
\end{equation}
with $\mu,\gamma,\delta>0$ and $\xi$ is a mean zero, stationary, Gaussian white noise process. 

The dynamics of periodically driven deterministic oscillators (i.e. $\sigma=0$ in \eqref{E:Cartmell}) has been studied extensively in the literature and is well understood: see for example, \cite{Bol64, Nay79} for \emph{weakly} nonlinear (i.e. $\gamma$ is of the order of $\eps \ll 1$  in \eqref{E:Cartmell}) systems. Deterministic oscillators with \emph{strong} nonlinearities and weak damping ($\delta \sim O(\eps)$) are studied as weakly perturbed Hamiltonian systems in, for example, \cite{Guc83} and \cite{Morozov1998}. 

On the other hand, in the presence of noise ($\sigma \neq 0$) and absence of periodic perturbations ($\alpha=0,\eta=0$), \eqref{E:Cartmell} represents the {\it noisy Duffing} equation. It has been studied with $\delta \sim O(\eps)$, $\sigma \sim O(\sqrt{\eps}))$ in \cite{Bol84, Sri:Ebe86, ArnL96, Lia99, SriSowVed}, to name a few. Let the `energy' of \eqref{E:Cartmell} be defined as $H_t:=\frac12\dot{q}_t^2-\frac{\mu}{2}q_t^2+\frac{\gamma}{4}q_t^4$. Employing the general technique of stochastic averaging developed in \cite{Sri:Fre98}, \cite{SriSowVed} showed that the distribution of the rescaled process $H^\eps_t:=H_{t/\eps}$ for $t\in [0,T]$ converges as $\eps \to 0$ to the distribution of a specific diffusion on a graph. For small $\eps$, this graph-valued process gives a good approximation to the dynamics of \eqref{E:Cartmell} on time intervals of order $1/\eps$.

Reference~\cite{Sri2001b} developed a unified approach for studying the dynamics of weakly nonlinear oscillators under \emph{both} periodic and noisy perturbations. It considers the case where $\mu<0$ (single well potential), and $\gamma, \eta, \sigma \sim O(\eps)$ and $\delta \sim O(\eps^2)$. In this case, \cite{Sri2001b} identified certain quantity\footnote{The system obtained by averaging out periodic oscillations due to $\eta$ has a constant of motion $H$. More precisely, let $\omg=\sqrt{-\mu}$ and  let $x^{(i)}_t$ be defined by $q_t=x^{(1)}_t\cos(\omg t)+x^{(2)}_t\sin(\omg t)$, $\dot{q}_t=-\omg x^{(1)}_t\sin(\omg t)+\omg x^{(2)}_t\cos(\omg t)$; then $H_t=\frac14((x^{(2)}_t)^2-(x^{(1)}_t)^2)+\frac{3\gamma}{16}((x^{(1)}_t)^2+(x^{(1)}_t)^2)^2$.} $H_t$ and, using non-standard averaging techniques, showed that the distribution of the rescaled process $H^\eps_t:=H_{t/\eps^2}$ for $t\in [0,T]$ converges as $\eps \to 0$ to the distribution of a specific diffusion on a graph. For small $\eps$, this graph-valued process gives a good approximation to the dynamics of \eqref{E:Cartmell} on time intervals of order $1/\eps^2$.
 
In this paper we study \eqref{E:Cartmell} with parameters scaled as:
\begin{equation}
\label{E:Cartmell_scale}
\ddot q_t - \mu  q_t +\gamma q_t^3=  \eps(\alpha \cos(\nu t)+ \mu\eta \cos(\nu t)q_t - \delta \dot{q}_t) + \eps^{\kappa}\sigma \,\xi(t),
\end{equation}
with $\mu,\gamma,\delta>0$ and $\xi$ is a mean zero, stationary, Gaussian white noise process and $\kappa\geq 1$.  
Equation \eqref{E:Cartmell_scale} can be studied as a perturbation of the Hamiltonian system 
\begin{equation}\label{eq:unpertsysHamsys}
\dot{q}=\frpt{H}{p},  \qquad \qquad \dot{p}=-\frpt{H}{q}
\end{equation} 
with the Hamiltonian 
\begin{equation}\label{eq:HamDefprimary}
H(q,p) = \frac12p^2 + U(q), 
\end{equation}
where $U$ is a double-well potential
\begin{equation}\label{eq:PotDefprimary}
U(q)=- \frac{\mu}{2} q^2  + \frac{\gamma}{4} q^4 .
\end{equation}

\begin{figure}\label{fig:potandcontourplot}
\centering
\includegraphics[scale=0.5]{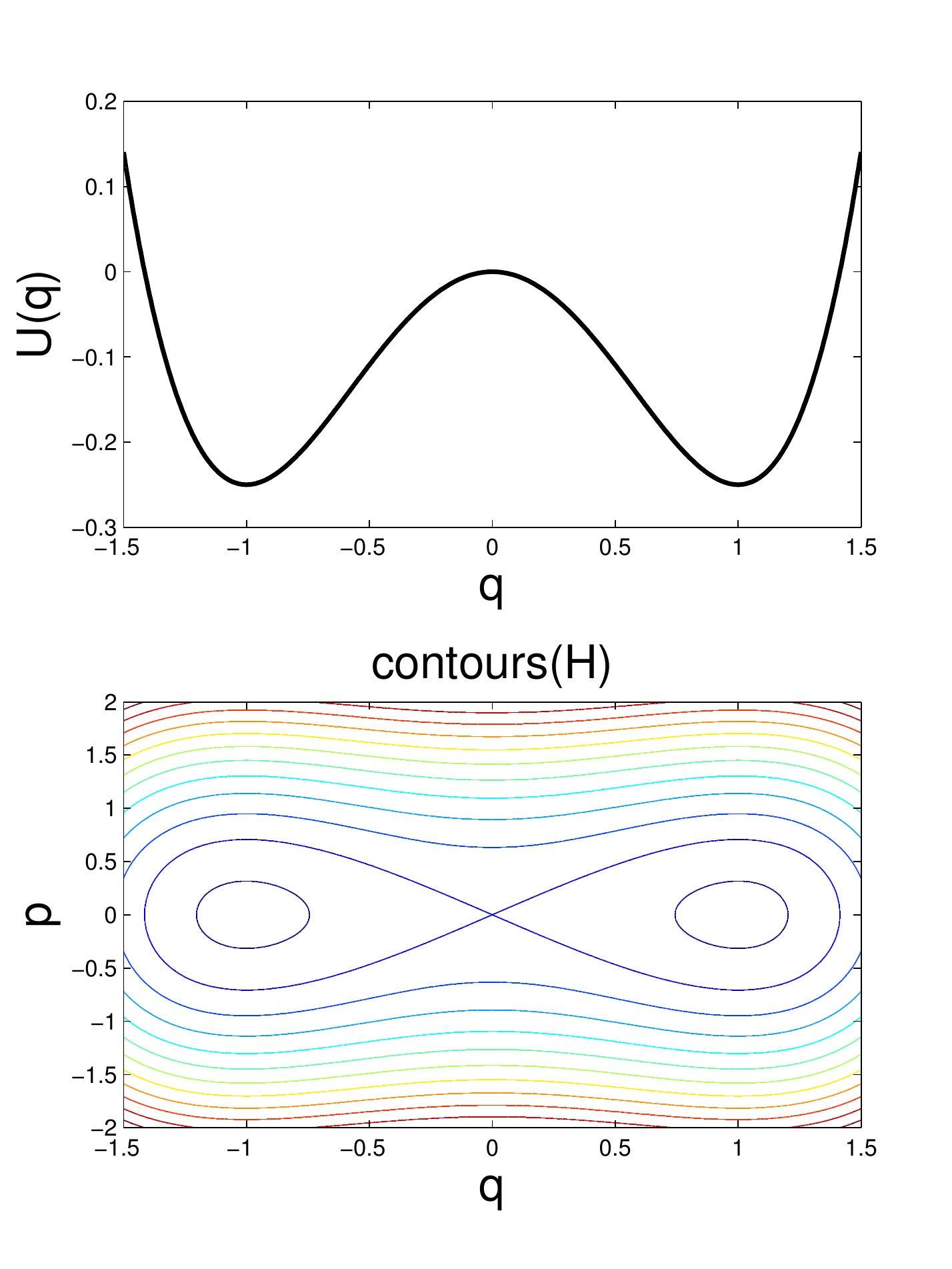}
\caption{Double-well potential $U$ and the contour plot for $H$.}
\end{figure}

The unperturbed system \eqref{eq:unpertsysHamsys}, written explicity as
\begin{equation}\label{eq:upexplicitly}
\ddot q_t - \mu  q_t +\gamma q_t^3=0,
\end{equation}
has three fixed points: the fixed point $q=0$ is a saddle and the other two fixed points $q=\pm\sqrt{\mu/\gamma}$ (corresponding to the bottom of the wells in the double-well potential) are centers. 
The shape of potential $U$ and the contour plot of $H$ is shown in the figure \ref{fig:potandcontourplot}. Note that $p$ is same as $\dot{q}$. We also use the notation $(q_1,q_2)$ for $(q,\dot{q})$.

Let $(I,\varphi)$ be action angle variables corresponding to the unperturbed system \eqref{eq:unpertsysHamsys} and assume the transformations
$$I=I(q_1,q_2), \qquad \varphi=\varphi(q_1,q_2),$$
$$q_1=q_1(I,\varphi), \qquad q_2=q_2(I,\varphi)$$ can be written. 
Then, the system \eqref{E:Cartmell_scale} with $\eps=0$ can be written as
\begin{equation}
\dot{I}=0, \qquad \qquad \dot{\varphi}=\Omg(I).
\end{equation}

Let $\theta$ be the angle variable corresponding to $\dot{\theta}=\nu$,  $\nu$ being the frequency of periodic excitations.
In the perturbed system~\eqref{E:Cartmell_scale} with $\eps \neq 0$ and $\sigma=0$, the frequency $\Omg(I)$ changes with time and if the frequencies $\Omg(I)$ and $\nu$ are non-commensurable, then the $(\varphi,\theta)$ orbits densely fill the state-space and the motion is called {\em quasi-periodic}. Resonance occurs when the frequencies $\nu$ and $\Omg(I)$ are commensurable or nearly commensurable and in this case orbits do not densely fill the state-space. Since $\Omg$ depends on the action $I$, the resonance will depend on certain values of the action. The region of state space which is close to the points where $\Omg(I)$ and $\nu$ are commensurable with a specific ratio is called a resonance zone.
The trajectories starting in some small set of initial conditions get `captured' into a resonance zone and those starting in other initial conditions `pass-through'. For detailed description of such phenomenon, see chapter 5 of \cite{ArnoldDynsys}. Clear understanding of the terms `resonance zone' and `capture' would be obtained by the end of section \ref{sec:captureIntoRes}.

Typically, studies which treat nonlinear oscillators as perturbation of Hamiltonian systems involve some kind of averaging principle. Issues that arise in obtaining averaging principle in presence of resonances are discussed, for example, in \cite{Morozov1998} and \cite{ArnoldDynsys}. In \cite{ArnoldDynsys} one can find a discussion on existing results about
\begin{itemize}
\item measure of the set of initial conditions that get captured into a resonance zone,
\item bounds on the error in approximating $I$ by an averaged system.
\end{itemize}
In studying the dynamics close to a resonance zone, partial averaging is employed as discussed in \cite{Morozov1998}. For example let $I_r$ be such that  $n\nu=m\Omg(I_r)$  where $n,m$ are integers. Then, \emph{the region of the state-space where $I$ is close to $I_r$ is called $m:n$ resonance zone}. In this region, introduce a new variable $\psi=\varphi-\frac{n}{m}\theta$. The dynamics in this region can be described using $(I,\psi)$ which are slow-variables while averaging out the fast variable $\theta$.

Periodically driven strongly nonlinear systems with noise are considered in \cite{FWMult}. It assumes that the noise in $(I,\varphi)$ variables is uniformly non-degenerate and obtains an averaging principle to the effect that the resonances could be totally ignored. However, the system \eqref{E:Cartmell_scale} that is being considered here, does not obey that hypothesis. In fact, any Hamiltonian system where the noise is added only to velocity coordinates does not obey that hypothesis---see the paragraph immediately following theorem 2.1 of \cite{FWMult} showing the restrictive nature of that hypothesis. Further, the noise considered in \cite{FWMult} is stronger than what we consider here: for the purpose of comparison, we take $\kappa\geq 1$ in \eqref{E:Cartmell_scale}, whereas the strength of noise in \cite{FWMult} corresponds to taking $\kappa=1/2$.

In light of the above discussion, the intention of this paper is to study the effect of weak noise on the escape from a resonance zone in \eqref{E:Cartmell_scale}. The phase space for the corresponding deterministic system (i.e. \eqref{E:Cartmell_scale} with $\sigma=0$) consists of many resonance zones in which some trajectories of the deterministic system can get `trapped'. When $\kappa>1$ the noise is very weak and so large deviations from the corresponding deterministic system  occur with very low probability.  The rate at which noise facilitates the `escape' from resonance is the subject of this paper. 
The material in this paper is organized as follows. 

In section \ref{sec:dynClose2ResZone} we derive the evolution equations for action-angle coordinates in the presence of noise. Then we make a change of variables which amounts to zooming in to a resonance zone.  
In section \ref{sec:captureIntoRes} we show the state-space of the dynamics in the resonance zone and state the well known problem of capture into resonance. We identify a variable $\mathcal{H}^\eps$ whose value can be used to indicate capture. 
In section \ref{sec:convgToDetm} we study the behaviour of $\mathcal{H}^\eps$ as $\eps\to 0$ and using averaging techniques identify the limit of $\mathcal{H}^\eps$ in the sense of distribution. For small $\eps$, the limiting process can be used to approximate the dynamics. In section \ref{sec:EscpFrmTrap} we find the mean time of exit from a resonance zone and study the dependence of mean-exit-time on the oscillator parameters.  It will be shown that the trajectories of the oscillator \eqref{E:Cartmell_scale} trickle down close to the bottom of the wells of the potential $U$. In section \ref{sec:resatbottom} we focus on the case where resonance occurs at the bottom of the potential wells. In such a scenario multiple solutions exist and we study how the rates of switching between the domains-of-attraction of different solutions depend on the oscillator parameters.

%##################
%##################
%##################
\section{Dynamics close to a resonance zone: Capture into resonance.}\label{sec:dynClose2ResZone}
In this section we derive stochastic evolution equations for action-angle coordinates and localize them near a resonance zone. Then we study the deterministic dynamics ($\sigma=0$) in the resonance zone and explain the capture phenomenon. In section \ref{sec:mainkappaGT1} we study how noise facilitates the escape from resonance zone.

We fix $\mu=1$ and $\gamma=1$ in \eqref{E:Cartmell_scale}; however, the ideas in this paper would be valid for any $\mu,\gamma>0$.

Rewriting \eqref{E:Cartmell_scale} in state-space form, we have
\begin{align}\label{eq:q1q2sde}
dq_{1,t}^{\eps}&=q_{2,t}^\eps dt,\\ \nonumber
dq_{2,t}^{\eps}&=( q_{1,t}^\eps-(q_{1,t}^\eps)^3) dt + \eps\left ( \eta \cos(\nu t)  q_{1,t}^{\eps}  + \alpha \cos(\nu t)  
- \delta q_{2,t}^{\eps} \right)dt + {\eps}^{\kappa}\sigma dW_{t},
\end{align}
where $W$ is a Wiener process. 
The transformation to the action-angle variables
$$I=I(q_1,q_2), \qquad \varphi=\varphi(q_1,q_2),$$
$$q_1=q_1(I,\varphi), \qquad q_2=q_2(I,\varphi)$$ 
can be written explicitly for each of the regions separated by the homoclinic orbit, see \ref{appsec:computeIJ}.
The system \eqref{eq:q1q2sde} with $\eps=0$ takes the form
\begin{equation}
\dot{I}=0, \qquad \qquad \dot{\varphi}=\Omg(I).
\end{equation}
When $\eps>0$, apart from the above frequency $\Omg(I)$, system \eqref{eq:q1q2sde} also depends on the frequency $\nu$ of the periodic excitation.
For arbitrary positive integers $m,n$ let $I_r$ be the value of action such that $m \Omg(I_r)=n\nu$. The region of state space where the action $I$ is close to $I_r$ is called $m:n$ resonance zone.
We study the dynamics of the system \eqref{eq:q1q2sde} in the region where $I$ is close to the resonant value $I_r$.

\begin{remark}
Here $r$ is short for resonance $m:n$. For notational convenience we use $\Omg_r=\Omg(I_r)$ and $\Omg'_r=\frpt{\Omg}{I}\big|_{I=I_r}$ etc.
\end{remark}

We use the It\^{o} formula to study how $I$ and $\varphi$ evolve for \eqref{eq:q1q2sde}. For this purpose, let $I_t^\eps:=I(q_{1,t}^{\eps},q_{2,t}^{\eps})$, $\varphi_t^\eps :=\varphi(q_{1,t}^{\eps},q_{2,t}^{\eps})$. Let the angle $\theta$ evolve according to $d\theta_t=\nu dt$ and define the angle
$$\psi_t^\eps := \varphi_t^\eps - \frac{n}{m}\theta_t.$$
Define 
$$g_2(u,v,\theta)\deff \eta \cos(\theta) u  + \alpha \cos(\theta)  
- \delta v, $$
$$\mfF(I,\varphi,\theta) \deff \frac{\partial I(q_1,q_2)}{\partial q_2}g_2(q_1,q_2,\theta)\bigg|_{q_1(I,\varphi),q_2(I,\varphi)},$$ 
$$\mfG(I,\varphi,\theta) \deff \frac{\partial \varphi(q_1,q_2)}{\partial q_2}g_2(q_1,q_2,\theta)\bigg|_{q_1(I,\varphi),q_2(I,\varphi)}.$$
Then, using the It\^{o} formula we get
\begin{align}
\label{eq:Ipsithe_stch}
\begin{cases}dI_t^\eps\,\,=\,\,\eps \mfF(I_t^\eps,\psi_t^\eps+\frac{n}{m}\theta_t,\theta_t)dt\,\, +\,\, \eps^\kappa\sigma \frpt{I}{q_2}\bigg|_{(I_t^\eps,\,\psi_t^\eps +\frac{n}{m}\theta_t)}dW_{t}  \,\,+\,\,\frac12\eps^{2\kappa}\sigma^2\frptsec{I}{q_2}\bigg|_{(I_t^\eps,\psi_t^\eps+\frac{n}{m}\theta_t)}dt, 
\\
d\psi_t^\eps\,\,=\,\, (\Omega(I_t^\eps)-\Omega_r)dt \,\,+\,\, \eps \mfG(I_t^\eps,\psi_t^\eps+\frac{n}{m}\theta_t,\theta_t)dt \,\,+\,\, \eps^\kappa\sigma \frpt{\varphi}{q_2}\bigg|_{(I_t^\eps,\psi_t^\eps+\frac{n}{m}\theta_t)}dW_{t} \\  \qquad \qquad \qquad  +\,\,\frac12\eps^{2\kappa}\sigma^2\frptsec{\varphi}{q_2}\bigg|_{(I_t^\eps,\psi_t^\eps+\frac{n}{m}\theta_t)}dt,   
\\
d{\theta}_t\,\,=\,\,\nu dt. 
\end{cases}
\end{align}
Here $\frpt{I}{q_2}\bigg|_{(I_t^\eps,\,\psi_t^\eps +\frac{n}{m}\theta_t)}$ denotes the partial derivative of the function $I(q_1,q_2)$ with respect to the second argument evaluated at $(q_1(I_t^\eps,\,\psi_t^\eps +\frac{n}{m}\theta_t),q_2(I_t^\eps,\,\psi_t^\eps +\frac{n}{m}\theta_t))$.

When $I_t^\eps$ is close to $I_r$, the difference $\Omega(I_t^\eps)-\Omega_r$ is small and hence $\psi_t^\eps$ evolves slowly compared to $\theta_t$. So, $(I,\psi)$ are slow variables and $\theta$ is a fast variable.

Since we are interested in the dynamics close to the resonance $I=I_r$, we make a change of variables in order to derive simpler equations that describe the dynamics in the resonance zone. Substituting the following standard space and time scaling (see, e.g. \cite{Neishtadt1991})
\begin{equation}\label{eq:rescalingdefinitions_hpsitheta0}
h_t^\eps\deff \frac{1}{\sqrt{\eps}}(I_{t/{\sqrt{\eps}}}^\eps-I_r), \qquad \psih_t^\eps\deff \psi_{t/{\sqrt{\eps}}}^\eps, \qquad \theta_t^\eps\deff \theta_{t/{\sqrt{\eps}}},
\end{equation}
into \eqref{eq:Ipsithe_stch} and Taylor-expanding in powers of $\sqrt{\eps}$ about $I_r$, we get, with higher order terms subsumed in $\mathfrak{R}$
\begin{align}\label{eq:hpsitheta_proprescale_a0}
dh_t^\eps &= \mfF dt + {\sqrt{\eps}}\, \mfF^{\prime}h_t^\eps dt +\eps^{\kappa-\frac{3}{4}} \,\sigma \frpt{I}{q_2}dW_{t}  + \mathfrak{R}_{1,t}^\eps dt +  \hat{\mathfrak{R}}_{1,t}^\eps dW_{t}, \\ \label{eq:hpsitheta_proprescale_b0}
d\psih_t^\eps &= \Omg^{\prime}_rh_t^\eps dt + {\sqrt{\eps}}\, \bigg(\Omg^{\prime\prime}_r\frac12(h_t^\eps)^2 + \mfG\bigg)dt  + \mathfrak{R}_{2,t}^\eps dt + \hat{\mathfrak{R}}_{2,t}^\eps dW_{t}, \\ \label{eq:hpsitheta_proprescale_c0}
d\theta_t^\eps &=  \frac{1}{\sqrt{\eps}}\nu dt, 
\end{align}
where $\prime$ indicates differentiation w.r.t $I$ and all terms (except $\mathfrak{R}$) are evaluated at $(I_r,\psih_t^\eps+\frac{n}{m}\theta_t^\eps,\theta_t^\eps).$ When $\kappa \geq 1$, the higher order terms are $\mathfrak{R}_{i}^\eps \sim O({\eps})$ and $\hat{\mathfrak{R}}_{i}^\eps\sim O(\eps^{\kappa-\frac{1}{4}})$, for $i=1,2$.

In the system \eqref{eq:hpsitheta_proprescale_a0}--\eqref{eq:hpsitheta_proprescale_c0}, $h_t^\eps,\psih_t^\eps$ are slow variables and $\theta_t^\eps$ is a fast variable.

%##################
%##################
%##################
\subsection{Capture into resonance}\label{sec:captureIntoRes}
From \eqref{eq:hpsitheta_proprescale_a0}-\eqref{eq:hpsitheta_proprescale_c0} it is clear that $\theta_t^\eps$ evolves at a faster rate than $h_t^\eps$ and $\psih_t^\eps$. In this section we show that, in the absence of noise ($\sigma=0$), averaging the fast $\theta$ oscillations would result in a Hamiltonian structure for $(\psi,h)$. Using the corresponding Hamiltonian $\mathcal{H}$, we explain the capture phenomenon.

For the purpose of averaging the fast $\theta$ oscillations, define an averaging operator $\la \cdot \ra$ as follows: for a function $f$ periodic in $\theta$ with period $2m\pi$ we define $\la f\ra = \frac{1}{2m\pi}\int_0^{2m\pi}f(\theta)d\theta$. Note that the functions $\theta \mapsto \mfF(I_r,\psi+\frac{n}{m}\theta,\theta)$ and $\theta \mapsto \mfG(I_r,\psi+\frac{n}{m}\theta,\theta)$ are periodic in $\theta$ with period $2m\pi$. To clearly indicate the dependence of the corresponding averaged function on $\psi$, we denote the averaged functions by $\la \mfF(\psi)\ra$ and $\la \mfG(\psi)\ra$.

For the analysis in this section, we neglect the stochastic term. To this end,  in \eqref{eq:hpsitheta_proprescale_a0}-\eqref{eq:hpsitheta_proprescale_c0} we set $\sigma=0$, ignore higher order terms $\mathfrak{R}$ and perform averaging with respect to $\theta$. Then we get 
\begin{equation}\label{eq:Ipsithe_stch_sc_timetransf_fin_det}
\left(\begin{array}{c}dh \\ d\psi \end{array}\right)=\left(\begin{array}{c}\la \mfF(\psi)\ra + \sqrt{\eps} \la \mfF'(\psi) \ra h \\ \Omg'_r h + \sqrt{\eps}(\frac12\Omg''_rh^2+\la \mfG(\psi)\ra ) \end{array}\right)dt,  
\end{equation}

The general structure of the averaged terms is  
\begin{align}\label{eq:genstructFG_a}
\la \mfF(\psi)\ra &=-\delta I_r + J_r\sin(m\psi/n), \\ \label{eq:genstructFG_b}
\la \mfF'(\psi)\ra &=-\delta + J^{\prime}_r\sin(m\psi/n), \\ \label{eq:genstructFG_c}
\la \mfG(\psi)\ra &= \frac{n}{m}J^{\prime}_r\cos(m\psi/n), 
\end{align}
where $$J_r=\eta J_r^\eta+\alpha J_r^\alpha, \qquad \quad J^{\prime}_r=\eta (J_r^\eta)^{\prime}+\alpha (J_r^\alpha)^{\prime}.$$
The method to obtain the above equations \eqref{eq:genstructFG_a}-\eqref{eq:genstructFG_c} and the quantities $J_r^\eta$, $J_r^\alpha, (J_r^\eta)', (J_r^\alpha)'$ is discussed in~\ref{appsec:computeIJ}.

\begin{figure}
\centering
\includegraphics[scale=0.32]{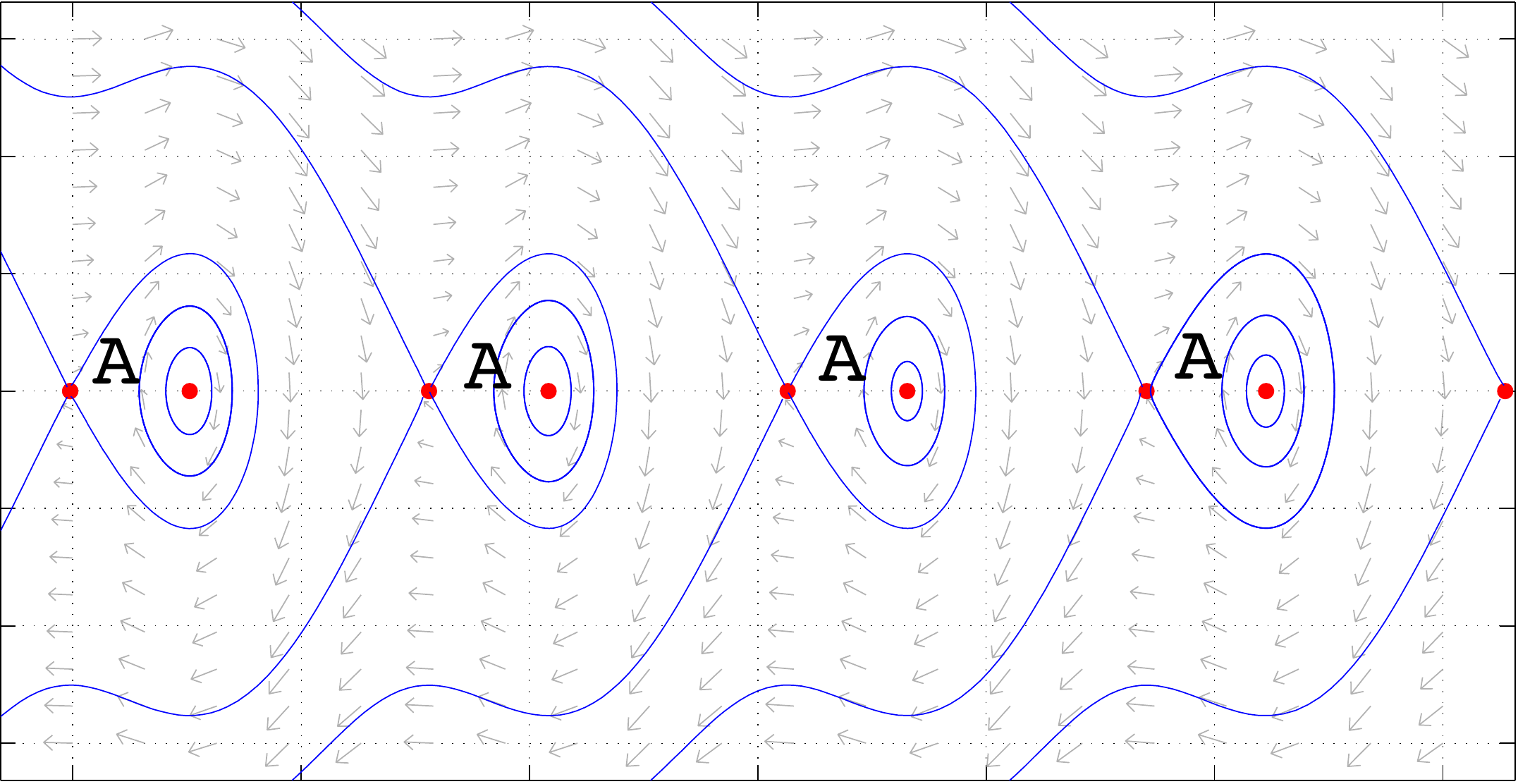}
\caption{Typical phase portrait for \eqref{eq:Ipsithe_stch_sc_timetransf_fin_det_basic} with $\Omg'_r>0$. Abscissa is $\psi$ and ordinate is $h$.}
\label{fig:phport_nodamp}
\end{figure}

\begin{figure}
\centering
\includegraphics[scale=0.35]{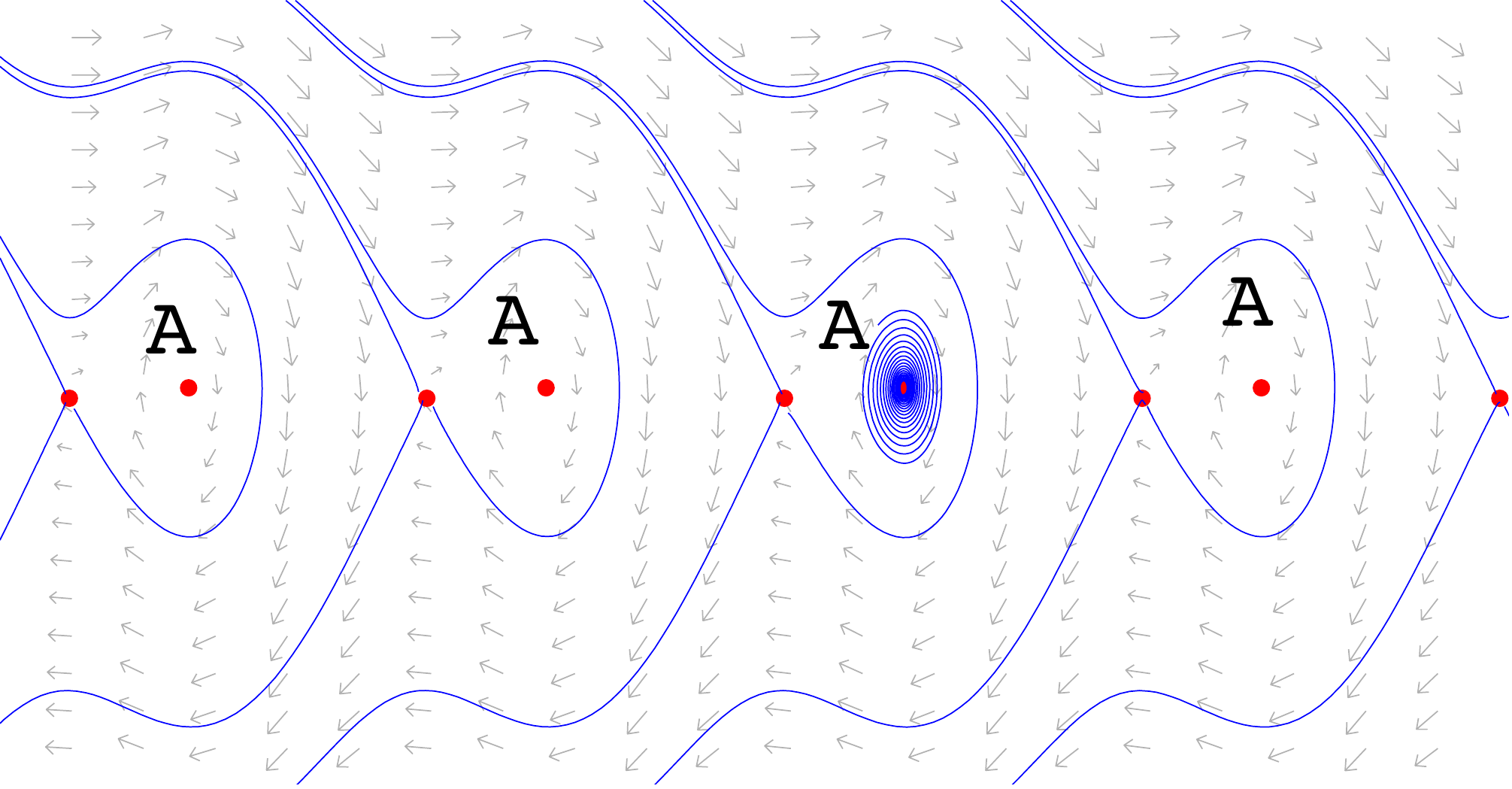}
\caption{Typical phase portrait for \eqref{eq:Ipsithe_stch_sc_timetransf_fin_det} with $\Omg'_r>0$. Abscissa is $\psi$ and ordinate is $h$. The system cannot leave the region $A$ in the absence of noise. The measure of the set of initial conditions that lead to trap in $A$ is small.}
\label{fig:phport2_damp}
\end{figure}

One of the reasons for localizing the equations near the resonance is that at $\eps=0$, equation~\eqref{eq:Ipsithe_stch_sc_timetransf_fin_det} reduces to a Hamiltonian system
\begin{equation}\label{eq:Ipsithe_stch_sc_timetransf_fin_det_basic}
\left(\begin{array}{c}dh \\ d\psi \end{array}\right)=\left(\begin{array}{c}\la \mfF(\psi)\ra  \\ \Omg'_r h \end{array}\right)dt,
\end{equation}
with the Hamiltonian
\begin{equation}\label{eq:HamOfAvgSysDef}
\mathcal{H}(\psi,h) = \frac12\Omg'h^2 - \int_{0}^{\psi} \la \mfF(\psi)\ra d\psi.
\end{equation}
Such Hamiltonians typically occur in resonant problems and \eqref{eq:HamOfAvgSysDef} represents a ``pendulum" under the action of an external torque~\cite{ArnoldDynsys, Neishtadt1991}. 
We can study \eqref{eq:Ipsithe_stch_sc_timetransf_fin_det} as a perturbation of \eqref{eq:Ipsithe_stch_sc_timetransf_fin_det_basic}.
Note that \eqref{eq:Ipsithe_stch_sc_timetransf_fin_det_basic} has a fixed point only if 
\begin{equation}\label{eq:nofixedpointif}
\delta I_r \leq |J_r|.
\end{equation}
The fixed points are given by the relation
$$\sin (m \psi/n) = \frac{\delta I_r}{J_r}, \qquad h=0.$$
There are many $\psi$ which satisfy the above equation.  Typical phase portraits (with $\Omg'_r>0$) for \eqref{eq:Ipsithe_stch_sc_timetransf_fin_det_basic}  and \eqref{eq:Ipsithe_stch_sc_timetransf_fin_det} are shown in the figures \ref{fig:phport_nodamp} and \ref{fig:phport2_damp}  respectively. 
The saddle(\emph{sd}) and center(\emph{sk}) fixed point pairs (i.e. the homoclinic orbit of the saddle encloses the center) for \eqref{eq:Ipsithe_stch_sc_timetransf_fin_det_basic} can be obtained as follows, with $j$ any integer and $\Psi_*=\frac{n}{m}\sin^{-1}(\delta I_r/J_r)$.
\begin{center}
  \begin{tabular}{ ||c || c || c ||}
\hline
\hline
	 & $\frac{m}{n}\Omg'_rJ_r\cos(m\Psi_*/n)>0$ &    $\frac{m}{n}\Omg'_rJ_r\cos(m\Psi_*/n)<0 $ \\
	\hline
	\hline
    $\Omg'_r>0$ & $\begin{array}{c}\psi_{sd}=\Psi_*+2j\pi\frac{n}{m}\\\psi_{sk}=\frac{n}{m}\pi-\Psi_*+2j\pi\frac{n}{m} \end{array}$ & $\begin{array}{c}\psi_{sk}=\Psi_*+2j\pi\frac{n}{m}\\\psi_{sd}=-\frac{n}{m}\pi-\Psi_*+2j\pi\frac{n}{m} \end{array}$ \\ \hline
    $\Omg'_r<0$ & $\begin{array}{c}\psi_{sd}=\Psi_*+2j\pi\frac{n}{m}\\\psi_{sk}=-\frac{n}{m}\pi-\Psi_*+2j\pi\frac{n}{m} \end{array}$ & $\begin{array}{c}\psi_{sk}=\Psi_*+2j\pi\frac{n}{m}\\\psi_{sd}=\frac{n}{m}\pi-\Psi_*+2j\pi\frac{n}{m} \end{array}$ \\ \hline
    \hline
  \end{tabular}
\end{center}
Consider figure \ref{fig:phport_nodamp}. All the fixed points have $h=0$.  Recall the definitions \eqref{eq:rescalingdefinitions_hpsitheta0}. Note that $h=0$ means $I=I_r$, i.e. the system is exactly at resonance. The figures \ref{fig:phport_nodamp} and \ref{fig:phport2_damp} show a finite region around $h=0$. In terms of $I$-coordinates, this region is a neighborhood of $I_r$ of a width of order $\sqrt{\eps}$. This is a \emph{resonance zone}.

A trajectory which starts at the top of the figure \ref{fig:phport2_damp} ($h>0$) but not in the narrow neck region would reach the bottom of the figure ($h<0$), i.e. the trajectory `passes' through the resonance zone. A trajectory which starts at the top of the figure \ref{fig:phport2_damp} ($h>0$)  in the narrow neck region enters the region $A$ and is trapped there. Let us call the region $A$ as `trap zone'.

For \eqref{eq:Ipsithe_stch_sc_timetransf_fin_det} the region marked $A$  (in figure \ref{fig:phport2_damp}) is a trap---the trajectories originating in $A$ cannot exit from it at all. However, for \eqref{eq:hpsitheta_proprescale_a0}--\eqref{eq:hpsitheta_proprescale_c0}, when $\sigma \neq 0$, the noise facilitates the escape.
We want to study how the noise facilitates the escape from the trap zone.

We denote by $\mathcal{H}|_{sd}$ the value of $\mathcal{H}$ evaluated at one saddle fixed point of \eqref{eq:Ipsithe_stch_sc_timetransf_fin_det_basic} and denote by $\mathcal{H}|_{sk}$ the value of $\mathcal{H}$ evaluated at the corresponding center fixed point of \eqref{eq:Ipsithe_stch_sc_timetransf_fin_det_basic}.

%###########################
%###########################
%###########################
\section{Stochastic dynamics close to a resonance zone.}
\label{sec:mainkappaGT1}

In this section we return back to the noisy system \eqref{eq:hpsitheta_proprescale_a0}--\eqref{eq:hpsitheta_proprescale_c0} and argue that $\mathcal{H}_t^\eps:=\mathcal{H}(\psih_t^\eps,h_t^\eps)$ is a good indicator of capture. In section \ref{sec:convgToDetm} we identify the limit of $\mathcal{H}_t^\eps$ as $\eps\to 0$ and use it to approximate the mean time of exit from the trap zone.

In \eqref{eq:hpsitheta_proprescale_a0}--\eqref{eq:hpsitheta_proprescale_c0}, to see the fluctuations of $\mathcal{H}(\psih^\eps_t,h^\eps_t)$, we need to look on an even longer $O(1/\sqrt{\eps})$ time scale. Hence we redefine the process $h,\psih,\theta$ process as using the following space and time scaling 
\begin{equation}\label{eq:rescalingdefinitions_hpsitheta}
h_t^\eps=\frac{1}{\sqrt{\eps}}(I_{t/{\eps}}^\eps-I_r), \qquad \psih_t^\eps=\psi_{t/{\eps}}^\eps, \qquad \theta_t^\eps=\theta_{t/{\eps}}.
\end{equation}
After doing a Taylor expansion about $I_r$, we get, with higher order terms subsumed in $\mathfrak{R}$
\begin{align}\label{eq:hpsitheta_proprescale_a}
dh_t^\eps &= \frac{1}{\sqrt{\eps}}\mfF dt + \mfF^{\prime}h_t^\eps dt +\eps^{\kappa-1}\sigma \frpt{I}{q_2}dW_{t}  + \mathfrak{R}_{1,t}^\eps dt +  \hat{\mathfrak{R}}_{1,t}^\eps dW_{t}, \\ \label{eq:hpsitheta_proprescale_b}
d\psih_t^\eps &= \frac{1}{\sqrt{\eps}}\Omg^{\prime}_rh_t^\eps dt + \bigg(\Omg^{\prime\prime}_r\frac12(h_t^\eps)^2 + \mfG\bigg)dt  + \mathfrak{R}_{2,t}^\eps dt + \hat{\mathfrak{R}}_{2,t}^\eps dW_{t}, \\ \label{eq:hpsitheta_proprescale_c}
d\theta_t^\eps &= \frac{1}{{\eps}}\nu dt, 
\end{align}
where $\prime$ indicates differentiation w.r.t $I$ and all terms (except $\mathfrak{R}$) are evaluated at $(I_r,\psih_t^\eps+\frac{n}{m}\theta_t^\eps,\theta_t^\eps).$ When $\kappa \geq 1$, the higher order terms are $\mathfrak{R}_{i}^\eps \sim O(\sqrt{\eps})$ and $\hat{\mathfrak{R}}_{i}^\eps\sim O(\eps^{\kappa-1/2})$, for $i=1,2$.

After averaging over $\theta$, the system \eqref{eq:hpsitheta_proprescale_a}-\eqref{eq:hpsitheta_proprescale_c} can be seen as a perturbation of the  Hamiltonian system \eqref{eq:Ipsithe_stch_sc_timetransf_fin_det_basic}. 
Let $\mathcal{H}_t^\eps:=\mathcal{H}(\psih_t^\eps,h_t^\eps)$, where $\mathcal{H}$ is defined in \eqref{eq:HamOfAvgSysDef}. The evolution of $\mathcal{H}_t^\eps$ can be obtained by applying the It\^{o} formula as
\begin{align} \label{eq:hpsitheta_proprescale_d}
d\mathcal{H}_t^\eps &= \frac{1}{\sqrt{\eps}}\Omg'_rh_t^\eps(\mfF-\la\mfF\ra)dt + \bigg((\Omg'_r\mfF'-\la \mfF \ra \frac12\Omg''_r)(h_t^\eps)^2-\la \mfF\ra \mfG\bigg)dt \\ \nonumber
 &  \qquad \quad +\frac12\eps^{2(\kappa-1)}\sigma^2 \Omg'_r\left(\frpt{I}{q_2}\right)^2dt + \eps^{\kappa-1}\sigma \Omg'_r h_t^\eps \frpt{I}{q_2}dW_{t} + \mathfrak{R}_{3,t}^\eps dt + \hat{\mathfrak{R}}_{3,t}^\eps dW_{t},  
\end{align}
where arguments for $\mfF$, $\la \mfF\ra$, $\mfG$, $\frpt{I}{q_2}$ are suppressed; and $\mathfrak{R}$ are higher order terms.
Since $\la \mfF-\la\mfF\ra \ra =0$, $\mathcal{H}_t^\eps$ evolves even slowly compared to $(\psih_t^\eps,h_t^\eps)$. 

Since our goal is to study the escape from the region marked $A$ on figure \ref{fig:phport2_damp}, we set the initial conditions to \eqref{eq:hpsitheta_proprescale_a}-\eqref{eq:hpsitheta_proprescale_b} in this region. In terms of $\mathcal{H}_t^\eps$ this amounts to specifying that $\mathcal{H}_0^\eps$ lies in between$\mathcal{H}|_{sk}$ and $\mathcal{H}|_{sd}$ (it can be shown that $\mathcal{H}|_{sd}>\mathcal{H}|_{sk}$ if $\Omg'_r>0$ and $\mathcal{H}|_{sd}<\mathcal{H}|_{sk}$ if $\Omg'_r<0$). When the noise is absent, i.e. $\sigma=0$, $\mathcal{H}_t^\eps$ drifts towards interior of the trap zone, i.e. towards $\mathcal{H}|_{sk}$. When $\sigma\neq 0$ the noise facilitates the escape. A good indicator of whether escape occured is $\mathcal{H}_t^\eps \geq \mathcal{H}|_{sd}$ in the case $\Omg'_r>0$ (if $\Omg'_r<0$ then a good indicator is $\mathcal{H}_t^\eps \leq \mathcal{H}|_{sd}$).  It is however not quite accurate for the following reason: The region of $(\psi,h)$ for which $\mathcal{H}$ lies between $\mathcal{H}|_{sk}$ and $\mathcal{H}|_{sd}$ is exactly the region marked $A$ in figure \ref{fig:phport_nodamp}. But the fixed points in figure \ref{fig:phport2_damp} differ from those of figure \ref{fig:phport_nodamp} by order $\sqrt{\eps}$. Hence the boundary of region $A$ in figure \ref{fig:phport2_damp} differs by a small amount from the boundary in figure \ref{fig:phport_nodamp}. Further, $\mathcal{H}_t^\eps$ could be a bit greater than $\mathcal{H}|_{sd}$ and still be in the small neck region which still leads to capture. Let $\mathcal{H}_*$ be the value for which we can be sure that escape occured if $\mathcal{H}_t^\eps \geq \mathcal{H}_*$. Then $\mathcal{H}|_{sd}$ differs from $\mathcal{H}_*$ by a very small amount that goes to zero as $\eps \to 0$. Keeping these caveats in mind, we still study the rate at which $\mathcal{H}_t^\eps$ exceeds $\mathcal{H}|_{sd}$ in presence of noise. However such transition is extremely unlikely because of the smallness of the noise. 

In section \ref{sec:convgToDetm} we show that, when $\kappa>1$, $\{\mathcal{H}_t^\eps\}_{t\in [0,T]}$ process for finite time $T$ converges in law as $\eps\to 0$ to a deterministic process for which $\mathcal{H}|_{sk}$ is a fixed point. More refined asymptotics in section \ref{sec:EscpFrmTrap} yields the mean time of exit from a trap zone.

%##############
%##############
\subsection{Convergence of $\{\mathcal{H}_t^\eps\}_{t\in [0,T]}$}\label{sec:convgToDetm}
In equation \eqref{eq:hpsitheta_proprescale_d}, note that $\la \mfF-\la\mfF\ra \ra = 0$. Hence, from the system of equations \eqref{eq:hpsitheta_proprescale_a}--\eqref{eq:hpsitheta_proprescale_d} it can be seen that $\mathcal{H}_t^\eps$ evolves slowly compared to $(\psih_t^\eps, h_t^\eps)$ which in-turn evolves slowly compared to $\theta_t^\eps$. Thus, to study the evolution of $\mathcal{H}_t^\eps$, we can average out the fast oscillations of $\theta_t^\eps$ and  $(\psih_t^\eps, h_t^\eps)$. The operator $\la \cdot \ra$ for averaging $\theta$ was already introduced in section \ref{sec:captureIntoRes}. For the purpose of averaging oscillations of $(\psih, h)$ along the Hamiltonian orbits, define an averaging operator $\Aavg$ as follows: 
\begin{definition}\label{def:Aavgdef}
For a function $f$ of $(\psih,h)$, the averaged function $\Aavg[f]$ is given by
$$\Aavg[f](\mathcal{\mfh})=\frac{1}{\mfT(\mathcal{\mfh})}\int_0^{\mfT(\mathcal{\mfh})}f(\psih(t),h(t))dt$$
where $(\psih(t),h(t))$ is the solution of the Hamiltonian system $\dot{\psih}=\frpt{\mathcal{H}}{h}, \dot{h}=-\frpt{\mathcal{H}}{\psih}$ with $\mathcal{H}(\psih,h)=\mfh$ and $\mfT(\mfh)$ is the time-period of the solution. The $\mfh$ is restricted to be in between $\mathcal{H}|_{sk}$ and $\mathcal{H}|_{sd}$; outside these values the orbit of the Hamiltonian system is not closed and the time-period is not defined.
\end{definition}
Since $\mfh$ is restricted to be in between $\mathcal{H}|_{sk}$ and $\mathcal{H}|_{sd}$ we define a stopping time
\begin{equation}\label{eq:stoptepsdef}
\stopteps:=\inf \{t> 0\,: \mathcal{H}_t^\eps \textrm{ is not in between } \mathcal{H}|_{sk} \textrm{ and } \mathcal{H}|_{sd}\}.
\end{equation}
More precisely, if $\Omg'_r>0$ then $\stopteps:=\inf \{t> 0\,: \mathcal{H}_t^\eps \geq \mathcal{H}|_{sd}\}$ and if $\Omg'_r<0$ then $\stopteps:=\inf \{t> 0\,: \mathcal{H}_t^\eps \leq \mathcal{H}|_{sd}\}$.

Using standard averaging techniques we have the following result.

\begin{theorem}\label{eq:WkConvgThmMain}
Let $\mathcal{A}_1(\psih,\theta)$ be defined by
\begin{equation}\label{eq:WkConvgThmMainA1def}
\nu \mathcal{A}_1(\psih,\theta)=\int_0^\theta\left(\mfF(I_r,\psih+\frac{n}{m}\theta,\theta)-\la \mfF(\psih) \ra\right)d\theta.
\end{equation}
Define $\mfB(\mfh)=\mfB_1(\mfh)+\mfB_2(\mfh)$ where $$\mfB_1=-\Omg'_r\Aavg\left[\left\la \mathcal{A}_1\mfF+\Omg'_rh^2\frpt{\mathcal{A}_1}{\psih} \right\ra\right], \qquad \quad \mfB_2=\Aavg \left[\left\la(\Omg'_r\mfF'-\la \mfF \ra \frac12\Omg''_r)h^2-\la \mfF\ra \mfG\right\ra \right].$$ In evaluating the $\theta$-averages $\la\ra$ in $\mfB_i$, the functions $\mfF$ and $\mfG$ should be treated as the maps $\theta \mapsto \mfF(I_r,\psih+\frac{n}{m}\theta,\theta)$, $\theta \mapsto \mfG(I_r,\psih+\frac{n}{m}\theta,\theta)$. The functions $\mfB_1,\mfB_2$ are evaluated in lemma \ref{lem:simplifyBXi}.
Let $\Xi=\sigma^2(\Omg'_r)^2\Aavg\left[\left\la \left(h \frpt{I}{q_2}\right)^2\right\ra \right]$ and $\mfB_{\sigma}=\frac12\sigma^2\Omg'_r\Aavg\left[\left\la \left(\frpt{I}{q_2}\right)^2\right\ra \right]$. 
Let the initial conditions to \eqref{eq:hpsitheta_proprescale_a}-\eqref{eq:hpsitheta_proprescale_c} be such that $\mathcal{H}_0^\eps$ lies in between $\mathcal{H}|_{sk}$ and $\mathcal{H}|_{sd}$. Then,
\begin{enumerate}
\item When $\kappa>1$, $\{\mathcal{H}_t^\eps\}_{t\in[0,T\wedge \stopteps]}$ converges to the deterministic process $\frac{d}{dt}\mathcal{H}_t^0 = \mfB(\mathcal{H}_t^0)$. 
\item When $\kappa=1$, $\{\mathcal{H}_t^\eps\}_{t\in[0,T\wedge \stopteps]}$ converges in law to the diffusion $\mathcal{H}_{t\wedge\stopt}^0$ given by
\begin{equation}\label{eq:avgkeq1diffproc}
d\mathcal{H}_t^0 = \left(\mfB(\mathcal{H}_t^0) + \mfB_{\sigma}(\mathcal{H}_t^0)\right)dt + \sqrt{\Xi(\mathcal{H}_t^0)}dW_{t},
\end{equation}
where $$\stopt:=\inf \{t> 0\,: \mathcal{H}_t^0 \textrm{ is not in between } \mathcal{H}|_{sk} \textrm{ and } \mathcal{H}|_{sd}\}.$$
\end{enumerate}
\end{theorem} 
{\bf Proof.} We write the infinitesimal generator of $(\psih_t^\eps, h_t^\eps, \theta_t^\eps)$ neglecting the higher order terms as
$$L^\eps=\eps^{-1}L_{0}+\eps^{-1/2}L_{1}+L_2+\eps^{2\rho}L_{\rho},$$
where $\rho=\kappa-1$ and $L_0=\nu\frpt{}{\theta}$, $L_1=\mfF\frpt{}{h}+\Omg'_rh\frpt{}{\psi}$, $L_2=\mfF'h\frpt{}{h}+(\Omg''_r\frac12h^2+\mfG)\frpt{}{\psi}$ and $L_{\rho}=\frac12\sigma^2(\frpt{I}{q_2})^2\frptsec{}{h^2}$. For real-valued functions $f_0$ defined on $[\mathcal{H}_{sk},\mathcal{H}_{sd}]$ define the operator
\begin{equation}\label{eq:avrgGenDeff}
\mathscr{L}=\left(\mfB_1(\mfh)+\mfB_2(\mfh)\right)\frpt{}{\mfh}+ \eps^{2(\kappa-1)}\left(\mfB_{\sigma}(\mfh)\frpt{}{\mfh}+\frac12\Xi(\mfh)\frac{\partial^2}{\partial \mfh^2}\right).
\end{equation}
We use the perturbed test function approach. Given smooth function $f_0$ of $\mathcal{H}$ we construct a function $f^\eps$ as
$$f^{\eps}(\psi,h,\theta)=f_0(\mathcal{H}(\psi,h))+\sqrt{\eps}f_1(\psi,h,\theta)+\eps f_2(\psi,h,\theta).$$
By appropriately choosing $f_1$ and $f_2$ we show that $\lim_{\eps \to 0}|L^\eps f^\eps - \mathscr{L}f_0|=0$ and $\lim_{\eps \to 0}|f^\eps - f_0|=0$. By theorem 4.8.2 in \cite{EthierKurtz} we would have that the finite-dimensional distributions of $\mathcal{H}^\eps$ converges to those of $\mathcal{H}^0$. For rigorous implementation of this procedure in an analogous situation see \cite{Sri2001b} and \cite{namSow}.

Expanding $L^\eps f^\eps$ we have
\begin{align*}
\frac{1}{\sqrt{\eps}}&\left[\nu\frpt{f_1}{\theta}+\left(\mfF\frpt{f_0}{h}+\Omg'_rh\frpt{f_0}{\psi}\right)\right]\\
\quad + \,\,&\left[\nu\frpt{f_2}{\theta}+\left(\mfF\frpt{f_1}{h}+\Omg'_rh\frpt{f_1}{\psi}\right)+\left(\mfF'h\frpt{f_0}{h}+(\Omg''_r\frac12h^2+\mfG)\frpt{f_0}{\psi}\right)\right]\\
\quad + \,\,\eps^{2(\kappa-1)}&\left[\frac{\sigma^2}{2}\left(\frpt{I}{q_2}\right)^2\frac{\partial^2f_0}{\partial h^2}\right] \,\,+\,\, O(\sqrt{\eps}).
\end{align*}
In the above expression we make use of the following two relations:
$$\frpt{f_0}{h}= \frpt{\mathcal{H}}{h}f'_0=\Omg'_rhf'_0, \qquad \frpt{f_0}{\psi}=\frpt{\mathcal{H}}{\psi}f'_0=-\la \mfF \ra f'_0.$$

Let $f_1=-\mathcal{A}_1\frpt{f_0}{h}+g$ where $\mathcal{A}_1$ is specified at \eqref{eq:WkConvgThmMainA1def} and $g$ be a function independent of $\theta$ which will be chosen in a moment. Then the $O(\frac{1}{\sqrt{\eps}})$-part in $L^\eps f^\eps$ becomes zero.

Let $\Phi_1=\mathcal{A}_1\mfF+\Omg'_rh^2\frpt{\mathcal{A}_1}{\psi}$, $\Phi_2=(\Omg'_r\mfF'-\la \mfF \ra \frac12\Omg''_r)h^2-\la \mfF\ra \mfG$, $\Phi_3=\frac12\sigma^2\left(\frpt{I}{q_2}\right)^2$ and $\Psi=(\Phi_2-\Omg'_r\Phi_1+\eps^{2(\kappa-1)}\Phi_3\Omg'_r)f'_0 + \eps^{2(\kappa-1)}\Phi_3(\Omg'_rh)^2f''_0$. Using the above $f_1$, the $O(1)$ and $O(\eps^{2(\kappa-1)})$ terms in $L^\eps f^\eps$ can together be written as
\begin{eqnarray*}
\nu\frpt{f_2}{\theta}+(\mfF-\la \mfF\ra)\left(\frpt{g}{h}-\mathcal{A}_1(\Omg'_r h)^2f''_0\right)+\Psi + \left[\la\mfF\ra\frpt{g}{h}+\Omg'_rh\frpt{g}{\psi}\right]=0.
\end{eqnarray*}
Recognizing that $\frpt{\mathcal{H}}{h}=\Omg'_rh$ and $\frpt{\mathcal{H}}{\psi}=-\la \mfF \ra$ we choose $g(\psi,h)$ as the solution of
\begin{equation*}
\frpt{\mathcal{H}}{h}\frpt{g}{\psi}-\frpt{\mathcal{H}}{\psi}\frpt{g}{h} = -\bigg(\la \Psi \ra - \Aavg[\la \Psi\ra]\bigg).
\end{equation*}
Note that LHS of the above equation is the derivative of $g$ along the Hamiltonian flow $\dot{\psi}=\frpt{\mathcal{H}}{h}$, $\dot{h}=-\frpt{\mathcal{H}}{\psi}$, and the average of RHS along a Hamiltonian orbit is zero.
We also choose
\begin{equation*}
\nu f_2(\psi, h,\theta) \,=\, -\int_0^\theta (\mfF-\la \mfF\ra)\left(\frpt{g}{h}-\mathcal{A}_1(\Omg'_r h)^2f''_0\right)d\theta \,-\,\int_0^\theta (\Psi-\la \Psi \ra)d\theta.
\end{equation*}
With the above choice of $f_1$ and $f_2$ we have that $L^\eps f^\eps = \Aavg[\la \Psi\ra] + O(\sqrt{\eps})$. Because $\Aavg[\la \Psi\ra]$ is same as $\mathscr{L}f_0$ we get that $L^\eps f^\eps -\mathscr{L}f_0$ is of order $\sqrt{\eps}$, and by construction $f^\eps-f_0$ is of order $\sqrt{\eps}$. $\hfill \square$
%##############
%##############

%##################
%##################
\subsection{Evaluation of $\mfB$ and $\Xi$}
The expression for $\mfB$ and $\Xi$ in theorem \ref{eq:WkConvgThmMain} can be simplified as follows.
\begin{lemma}\label{lem:simplifyBXi}
Let $g$ be defined by $g(\mfh)=\Aavg[h^2](\mfh)$. Then the coefficients $\mfB$ and $\Xi$ in theorem \ref{eq:WkConvgThmMain} can be simplified as
\begin{eqnarray*}
\mfB_1(\mfh)=0, \qquad \qquad \mfB_2(\mfh)=-\delta \Omg'_rg(\mfh), \\
\mfB_{\sigma}(\mfh)=\frac12\sigma^2\frac{\Omg'_rI_r}{\Omg_r}, \qquad
\Xi(\mfh)=\frac{\sigma^2(\Omg'_r)^2 I_r}{\Omg_r}g(\mfh).
\end{eqnarray*}
\end{lemma}
{\bf Proof.}
The simplified expressions for $\mfB_{\sigma}$ and $\Xi$ follow from the equality: (i) $\frpt{I}{q_2}=\frpt{I}{H}\frpt{H}{q_2}=\frac1\Omg q_2$ and (ii) by the definition of action $I_r=\frac{1}{2\pi}\oint q_2dq_1 = \frac{1}{2\pi/T}\frac1T\int_0^T q_2^2dt = \frac{1}{\Omg_r}\la q_2^2\ra$

It turns out that $\mfB_1 \equiv 0$. We do not know any general reason why it should be zero. However we prove this by evaluating the terms. Since it is quite tedious we shift it to~\ref{appsec:showingB1eq0}.
 
As for $\mfB_2$, using the structure of $\mfF$ from \eqref{eq:genstructFG_a}-\eqref{eq:genstructFG_c}, we get that
\begin{eqnarray} \label{eq:B1eval_aux} 
\mfB_2=-\delta(\Omg'_r-\frac12\Omg''_rI_r)\Aavg[h^2] + (\Omg'_rJ'_r-\frac12\Omg''_rJ_r)\Aavg[h^2\sin(m\psih/n)] \\  \nonumber
 \qquad \qquad - J_r'\frac{n}{m} \Aavg[(-\delta I_r+J_r\sin(m\psi/n))\cos(m\psih /n)]. 
\end{eqnarray}
Recalling the definition \ref{def:Aavgdef} of the operator $\Aavg$ we have
\begin{align}\label{eq:intbypartAux1}  
\Aavg[h^2\sin(m\psi/n)]&=\frac{1}{\mfT}\int_0^\mfT h^2(t)\sin(m\psi(t)/n)dt \\ \nonumber
&=\frac{1}{\mfT}\int_0^\mfT h^2(t)\left(\frac{\frac{dh}{dt}+\delta I_r}{J_r}\right)dt \\  \nonumber
&= \frac{1}{J_r\mfT}\frac13(h^3(\mfT)-h^3(0)) + \frac{\delta I_r}{J_r}\Aavg[h^2] = \frac{\delta I_r}{J_r}\Aavg[h^2]. 
\end{align}
Similarly, using $-\delta I_r+J_r\sin(m\psi/n)=\frac{dh}{dt}$ and doing integration by parts while using $\frac{d\psi}{dt}=\Omg'_rh$ yields
\begin{equation}\label{eq:intbypartAux2}
\Aavg[(-\delta I_r+J_r\sin(m\psi/n))\cos(m\psi/n)]=\frac{m}{n}\Omg'_r(\delta I_r/J_r)\Aavg[h^2].
\end{equation}
Employing the above two results in \eqref{eq:B1eval_aux} gives $\mfB_2(\mfh)=-\delta \Omg'_rg(\mfh)$. $\hfill \square$

%##############
%##############
%##############

\subsection{Escape from the trap zone. Case $\kappa>1$.}\label{sec:EscpFrmTrap}

In this section we derive an approximation to mean time of exit from trap zone using the limit processes $\mathcal{H}^0$ of theorem \ref{eq:WkConvgThmMain}. Then we study the dependence of mean exit time on the oscillator parameters.

Since we are interested in the escape from the trap zone (region $A$ in the figure \ref{fig:phport2_damp}), we consider the mean exit time $\expt_{\hbarst}[\stopteps]$ where $\stopteps$ is defined in \eqref{eq:stoptepsdef} and $\hbarst$ indicates that the initial condition is such that $\mathcal{H}^\eps_0=\hbarst$. We restrict ourselves to the case that $\hbarst$ lies between $\mathcal{H}|_{sk}$ and $\mathcal{H}|_{sd}$. 

Theorem \ref{eq:WkConvgThmMain} shows that, for $\kappa=1$ the distribution of $\mathcal{H}^\eps$ converges to that of $\mathcal{H}^0$ given by the stochastic differential equation (SDE) \eqref{eq:avgkeq1diffproc}, and for $\kappa>1$ the distribution converges to that of the deterministic equation $\frac{d}{dt}\mathcal{H}_t^0 = \mfB(\mathcal{H}_t^0)$. For small $\eps$ the distribution of $\mathcal{H}^0$ can be used to approximate the dynamics of $\mathcal{H}^\eps$. Note that the averaged generator \eqref{eq:avrgGenDeff} corresponds to the SDE
\begin{equation}\label{eq:avgkgt1diffproc}
d\mathcal{H}_t = \left(\mfB(\mathcal{H}_t) + \eps^{2(\kappa-1)}\mfB_{\sigma}(\mathcal{H}_t)\right)dt + \eps^{(\kappa-1)}\sqrt{\Xi(\mathcal{H}_t)}dW_{t}.
\end{equation}
In the case $\kappa=1$, the SDE \eqref{eq:avgkgt1diffproc} coincides with \eqref{eq:avgkeq1diffproc}. In the case $\kappa>1$, for small $\eps$, we expect the SDE \eqref{eq:avgkgt1diffproc} to give a better approximation to the dynamics than the deterministic equation $\frac{d}{dt}\mathcal{H}_t^0 = \mfB(\mathcal{H}_t^0)$. We have the following result for exit times of \eqref{eq:avgkgt1diffproc}.
\begin{theorem}\label{thm:exitTimeResult}
Let $\tau^\eps$ be the exit time of $\mathcal{H}_t$ defined by \eqref{eq:avgkgt1diffproc} from the region bounded by $\mathcal{H}_{sk}$ and $\mathcal{H}_{sd}$. Then, for $\kappa>1$,
\begin{equation}\label{eq:exitTimeResult}
\lim_{\eps\to 0}\eps^{2(\kappa-1)}\log \expt_{\hbarst}[\tau^\eps]\,\, =\,\, \msV(\mathcal{H}|_{sd}),
\end{equation}
where 
\begin{equation*}
\msV(\mfh)=\frac{2\delta \Omg_r}{\sigma^2\Omg'_rI_r}(\mfh-\mathcal{H}|_{sk}).
\end{equation*}
\end{theorem}
{\bf Proof.} Let $g$ be as defined in lemma \ref{lem:simplifyBXi}. Linearizing the Hamiltonian dynamics near $\mathcal{H}_{sk}$ it can be shown that close to $\mathcal{H}_{sk}$ the function $g$ behaves as  $g(\mfh)\approx \frac{\mfh-\mathcal{H}_{sk}}{\Omg'_r}\left[1-\frac{m^2}{8n^2}\frac{\Omg'_r(\mfh-\mathcal{H}_{sk})}{|\Omg'_rJ_r\frac{m}{n} \cos(m\Psi_*/n)|}\right]$. Using this fact, standard calculations (see table 6.2 in chapter 15 of \cite{KarTay} or chapter 8 in \cite{EthierKurtz}) show that $\mathcal{H}_{sk}$ is an entrance boundary for \eqref{eq:avgkgt1diffproc}. Let $\mathcal{L}$ be as defined in  \eqref{eq:avrgGenDeff}. Let $u(\mfh)=\expt_{\mfh}[\tau^\eps]$. Then $u$ is the solution of (see chapter 15.3 of \cite{KarTay})
\begin{equation*}
\mathcal{L}u = -1, \qquad u(\mathcal{H}_{sd})=0, \qquad u(\mathcal{H}_{sk})<\infty.
\end{equation*}
The above equation can be solved as
\begin{equation*}
u(\mfh)\,=\,\frac{\lambda}{\eps^{2\rho}\,\delta \Omg'}\int_{\mfh}^{\mathcal{H}_{sd}}dy\,\,e^{\lambda y/\eps^{2\rho}}\int_{\mathcal{H}_{sk}}^y d\eta\,\,e^{-\lambda\eta/\eps^{2\rho}} \Theta(y,\eta),
\end{equation*}
where $\rho=2(\kappa-1)$, $\lambda=\frac{2\delta\Omg_r}{\sigma^2\Omg'_rI_r}$ and $\Theta(y,\eta)=\exp\left(-\frac{1}{\Omg'_r}\int_{\eta}^y\frac{d\xi}{g(\xi)}\right)\frac{1}{g(\eta)}$. Using the behaviour of $g$ close to $\mathcal{H}_{sk}$ it can be shown that $\lim_{\eta \to \mathcal{H}_{sk}}\Theta(y,\eta)=:\Theta(y, \mathcal{H}_{sk})$ is finite. 

In the limit $\eps\to 0$, using Laplace's principle \cite{Olver1997}, $u$ can be approximated as
\begin{align*}
u(\mfh)\,&\approx\,\frac{\lambda}{\eps^{2\rho}\,\delta \Omg'}\int_{\mfh}^{\mathcal{H}_{sd}}dy\,\,e^{\lambda y/\eps^{2\rho}}\Theta(y, \mathcal{H}_{sk})\frac{\eps^{2\rho}}{\lambda}e^{-\lambda\mathcal{H}_{sk}/\eps^{2\rho}}\\
&\approx\,\frac{\lambda}{\eps^{2\rho}\,\delta \Omg'}\Theta(\mathcal{H}_{sd}, \mathcal{H}_{sk})\frac{\eps^{2\rho}}{\lambda}e^{\lambda \mathcal{H}_{sd}/\eps^{2\rho}}\frac{\eps^{2\rho}}{\lambda}e^{-\lambda\mathcal{H}_{sk}/\eps^{2\rho}}\\
&=\eps^{2(\kappa-1)}\frac{\Theta(\mathcal{H}_{sd}, \mathcal{H}_{sk})}{\lambda\delta\Omg'_r}\exp\left(\frac{\lambda(\mathcal{H}_{sd}-\mathcal{H}_{sk})}{\eps^{2(\kappa-1)}}\right)
\end{align*}
from which the desired result follows. $\hfill \square$

\vspace{40pt}
Note that $\expt[\tau^\eps]$ is of the order $e^{\msV(\mathcal{H}|_{sd})/\eps^{2(\kappa-1)}}$. Though the averaging result of theorem \ref{eq:WkConvgThmMain} holds only on times of order 1, because the fixed point $\mathcal{H}_{sk}$ is stable, we might expect that $\expt[\stopteps]$ would be approximately of the same order as $\expt[\tau^\eps]$.

The quantity $\msV(\mathcal{H}|_{sd})$ gives a measure of difficulty of escape from the trap zone:
\begin{equation}
\msV(\mathcal{H}|_{sd})= \frac{2\delta \Omg_r}{\sigma^2\Omg'_rI_r}(\mathcal{H}|_{sd}-\mathcal{H}|_{sk}).
\end{equation}
More precisely, it can be evaluated to be
\begin{align*}
\msV(\mathcal{H}|_{sd}) &= \frac{2 \Omg_r(n/m)}{\sigma^2|\Omg'_r|}\delta^2\left(2\sin^{-1}|\chi| -\pi + 2\frac{\sqrt{1-|\chi|^2}}{|\chi|}\right), \\  \chi &:= \frac{\delta I_r}{J_r}=\frac{\delta I_r}{\eta J_\eta + \alpha J_\alpha}.
\end{align*}
The condition \eqref{eq:nofixedpointif} entails that $|\chi|<1$ for the resonance zone to exist.
 
Since the function in the brackets is monotonically decreasing in $|\chi|$, we can deduce that for a fixed $\delta$, $\msV(\mathcal{H}|_{sd})$ is monotonically increasing in $|J_r|$, i.e. the higher the strength of periodic excitations the more difficult is the escape from the trap. For a fixed $J_r$, $\msV(\mathcal{H}|_{sd})$ has a unique maximum as a function of $\delta$. As $\delta$ increases to $\frac{|J_r|}{I_r}$, $\msV(\mathcal{H}|_{sd})$ decreases to $0$, because the area of the trap zone decreases to zero. As $\delta$ decreases to $0$, $\msV(\mathcal{H}|_{sd})$ also decreases to zero---this behaviour is not intuitive. Hence, for a fixed strength of periodic excitations, both high and low damping make the escape easier---intermediate values of damping make the escape difficult.

%################ AFTER ESCAPE ###########
%#####################################
\subsection{Post escape from the trap. Case $\kappa>1$.}
Immediately outside of the trap region A, the deterministic dynamics alone is enough to take the system out of the resonance zone (see the vector field in figure \ref{fig:phport2_damp}). Since the noise is small, getting re-trapped is a rare event, i.e. the system moves out of the resonance zone quickly. Once outside of the resonance zone, full-averaging i.e. averaging w.r.t $(\varphi,\theta)$ can be done. The full-averaged system shows that damping results in decrease of $I$ with time. However as $I$ decreases the system might enter a different resonance zone---from results of \cite{ArnoldDynsys} we know that the measure of the set of initial conditions which get trapped is small. Those that get trapped, escape at a rate governed by the results of section \ref{sec:EscpFrmTrap}. In such fashion the system evolves until it reaches a vicinity of $(q_1,q_2)=(\pm\sqrt{\mu/\gamma},0)$, i.e. a bottom of the wells in the potential $U$ of \eqref{eq:PotDefprimary}.

Note that we have not analysed the behaviour near the homoclinic orbit. So, the description in the above paragraph is valid for those trajectories which start within the region bounded by the homoclinic orbit in figure \ref{fig:potandcontourplot}. However, the analysis in previous sections is valid also for the resonance zones that lie outside the region bound by the homoclinic orbit.

 If the action at the bottom of the well $I_b:=I|_{q_1=\sqrt{\mu/\gamma},q_2=0}$ is such that  $\Omg(I_b)$ is in resonance with $\nu$, then interesting dynamics occurs. Such a situation is discussed in \cite{Dykman} in an attempt to explain phase-flip of electrons in external fields. For the sake of completeness we discuss this in section \ref{sec:resatbottom} when $\nu \approx 2\Omg(I_b)$. 

%#####################################
%################ kappa=1 ###########
%#####################################
\subsection{On the possibility of obtaining a large deviations principle for $\mathcal{H}^\eps$ when $\kappa>1$.} \label{sec:possibleLDP}
The convergence of $\mathcal{H}^\eps$ to the deterministic process $\frac{d}{dt}\mathcal{H}_t^0 = \mfB(\mathcal{H}_t^0)$ and the exponential order of the mean time of exit (which is a large deviation from the deterministic system) raise the question whether a large deviations principle (LDP) can be obtained for $\mathcal{H}^\eps$ process. Note that the limit $\mathcal{H}^0$ was obtained by averaging. The interplay of averaging and large deviations is studied in, for example, \cite{SowersLDP}, \cite{VeretLDPAvg2000}, \cite{SpilioMultiscaleLDP} and \cite{SpilioHu}. Reference \cite{SpilioMultiscaleLDP} considers slow-fast systems with two time-scales whose diffusion coefficient is uniformly non-degenerate and shows that, under certain conditions, the slow system possesses an LDP with the rate function same as the rate function of the slow process obtained by averaging out the fast process (see the paragraph after theorem 3.4 in \cite{SpilioMultiscaleLDP}). The system of equations \eqref{eq:hpsitheta_proprescale_a}-\eqref{eq:hpsitheta_proprescale_d} is a three time-scale system, i.e. $\theta^\eps$ evolves at a faster rate than $(\psih^\eps,h^\eps)$ which in turn evolve at a faster rate than $\mathcal{H}^\eps$. Further, the diffusion coefficient in \eqref{eq:hpsitheta_proprescale_a}-\eqref{eq:hpsitheta_proprescale_d} is not non-degenerate. The diffusion coefficient for the averaged process, $\sqrt{\Xi}$ in \eqref{eq:avgkgt1diffproc}, behaves near the sink fixed point as $\sqrt{\Xi(\mfh)}\approx \sqrt{\Omg'(\mfh-\mathcal{H}|_{sk})}$, which also is not uniformly non-degenerate. Large deviations for processes with diffusion coefficients behaving as $\sigma(x)\approx \sqrt{x}$ are studied in \cite{FischerChiarini}. However \cite{FischerChiarini} assumes that the drift coefficient is positive when diffusion coefficient is zero. This does not hold for \eqref{eq:avgkgt1diffproc} because $\mfB(\mathcal{H}|_{sk})=0$. A preliminary analysis shows that a study of amount of time spent near the boundary $\mathcal{H}|_{sk}$ is crucial in order to obtain an LDP.

%#####################################
%################ kappa=1 ###########
%#####################################
\subsection{Diffusion limit: $\kappa=1$} \label{sec:averging}
When $\kappa=1$, for small $\eps$, the dynamics of $\mathcal{H}^\eps$ can be approximated by that of the diffusion process \eqref{eq:avgkeq1diffproc}. However, the dynamics in this case is difficult to understand because (i) even after an exit from the trap zone, the noise is strong enough to make the system re-enter the trap zone; (ii) averaging over $(\psih,h)$ cannot be done outside the trap zone.

%#####################################
%################ kappa=1 ###########
%#####################################
\section{Resonance at the bottom of the potential well}\label{sec:resatbottom}
As mentioned earlier, trajectories of the oscillator \eqref{eq:q1q2sde} trickle down close to the point $q=\sqrt{\mu/\gamma}$, i.e. the bottom of the potential well for the unperturbed system \eqref{eq:upexplicitly}. Linearizing \eqref{eq:upexplicitly} about $q=\sqrt{\mu/\gamma}$, i.e.  setting $x=q-\sqrt{\mu/\gamma}$ and retaining terms linear in $x$, we get $\ddot{x}+2\mu x=0$. This shows that, close to the bottom of the potential well, the unperturbed system behaves approximately like a harmonic oscillator with frequency $\sqrt{2\mu}$. Interesting dynamics occurs if the forcing frequency $\nu$ in the perturbed system \eqref{eq:q1q2sde} is in resonance with the oscillator natural frequency $\sqrt{2\mu}$. In the absence of noise ($\sigma=0$) multiple solutions co-exist for \eqref{eq:q1q2sde} with the state-space partitioned as domains of attractions for individual solutions. When $\sigma\neq 0$, the noise facilitates switching of the trajectories between the domains of attraction. 

In this section we study the perturbed system \eqref{eq:q1q2sde} in the case where the forcing frequency $\nu$ is close to $2\sqrt{2\mu}$, i.e. we investigate $2:1$ resonance.  We assume $\nu = 2\sqrt{2\mu}(1+\eps \lambda)$ where $\lambda$ is a detuning parameter. Such a situation is discussed in \cite{Dykman} in an attempt to explain phase-flip of electrons in external fields. We study the dependence of switching rates on oscillator parameters.

Using the transformation $x_{1,t}^\eps=q_{1,t}^{\eps}-\sqrt{\mu/\gamma}$ and $x_{2,t}^\eps=q_{2,t}^{\eps}$ in equation \eqref{eq:q1q2sde} we find
\begin{equation}\label{eq:q1q2sde_auxX}
\hspace{-40pt}\begin{cases}dx_{1,t}^{\eps}=x_{2,t}^\eps dt,\\
dx_{2,t}^{\eps}=-(2\mu x_{1,t}^\eps+3\gamma \sqrt{\mu/\gamma}(x_{1,t}^\eps)^2+\gamma (x_{1,t}^\eps)^3) dt  + \eps \sqrt{\mu/\gamma}\eta \mu \cos(\nu t)dt\\
 \qquad \qquad + \eps\left (\eta \mu \cos(\nu t)  x_{1,t}^{\eps}  + \alpha \cos(\nu t)  
- \delta x_{2,t}^{\eps} \right)dt + {\eps}^{\kappa}\sigma dW_{t}.\end{cases}
\end{equation}
The forcing $\eps \sqrt{\mu/\gamma}\eta \mu \cos(\nu t)$ induces an approximately periodic motion with an amplitude of the order $O(\eps)$. However, significant length scale in the system turns out to be of the order $O(\sqrt{\eps})$. Further, the system \eqref{eq:q1q2sde_auxX} can be simplified by performing a near identity transformation which eliminates the quadratic nonlinearities in \eqref{eq:q1q2sde_auxX}. This motivates the following sequence of transformations on \eqref{eq:q1q2sde_auxX}:
$$\begin{cases}v_{1,t}^{\eps}=\frac{1}{\sqrt{\eps}}\left(x_{1,t}^\eps+\eps \frac{\etat\mu \sqrt{\mu/\gamma}}{\nu^2-2\mu}\cos(\nu t)\right), \\
v_{2,t}^{\eps}=\frac{1}{\sqrt{\eps}}\frac{1}{\sqrt{2\mu}}\bigg(x_{2,t}^{\eps}-\eps\frac{\etat\mu \sqrt{\mu/\gamma}}{\nu^2-2\mu}\nu\sin(\nu t) \bigg),\end{cases}$$
$$\left(\begin{array}{c}y_{1,t}^{\eps}\\y_{1,t}^{\eps}\end{array}\right)=\left(\begin{array}{c}v_{1,t}^{\eps}\\v_{2,t}^{\eps}\end{array}\right)+\sqrt{\eps}\sqrt{\frac{\gamma}{4\mu}}\left(\begin{array}{c}(v_{1,t}^{\eps})^2+2(v_{2,t}^{\eps})^2\\-2v_{1,t}^{\eps}v_{2,t}^{\eps}\end{array}\right)+{\eps}{\frac{\gamma}{4\mu}}\left(\begin{array}{c}(2v_{1,t}^{\eps})^3-3v_{1,t}^{\eps}(v_{2,t}^{\eps})^2\\3(v_{1,t}^{\eps})^2v_{2,t}^{\eps}\end{array}\right).$$
We find that the dominant dynamics of $y$ is rotation with frequency close to $\frac12\nu$. So, we make one additional transformation to remove the rotation:
\begin{equation*}
\left(\begin{array}{c}z_{1,t}^{\eps}\\z_{1,t}^{\eps}\end{array}\right)=e^{-tB/\eps}\left(\begin{array}{c}y_{1,t/\eps}^{\eps}\\y_{2,t/\eps}^{\eps}\end{array}\right), \qquad \qquad B=\frac12\nu\left(\begin{array}{cc}0 & 1 \\ -1 & 0\end{array}\right).
\end{equation*}
Then we get, 
\begin{eqnarray}\nonumber
\hspace{-50pt}dz_{t}^{\eps}=e^{-tB/\eps}\left\{\left(\begin{array}{cc}0 & 0\\(\etat \sqrt{2\mu})\cos(\nu t/\eps)& -\deltat \end{array}\right)-\frac{3\gamma}{4\mu} (z_1^2+z_2^2)\frac{2}{\nu}B -\lambda B\right\}e^{tB/\eps}z_{t}^{\eps}dt \\ 
 \qquad \qquad + \eps^{\kappa-1}\frac{\sigmat}{\sqrt{2\mu}}e^{-tB/\eps}\left(\begin{array}{c}0\\ 1\end{array}\right)dW_{t}\qquad + \,\,h.o.t \label{eq:RDS}
\end{eqnarray}
The higher order terms are not significant for dynamics of $z^\eps$ on a fixed time interval $[0,T]$. The fast oscillating coefficients in the above equation can be averaged out.
Define the averaged drift coefficient by
\begin{equation*}
\msB(z)=\left(-\frac{3\gamma}{4\mu}(z_1^2+z_2^2)-\frac{\nu}{2}\lambda \right)\left(\begin{array}{c}z_2\\ -z_1\end{array}\right)-\frac{\deltat}{2}\left(\begin{array}{c}z_1\\z_2\end{array}\right)+\frac{\etat \sqrt{2\mu}}{4}\left(\begin{array}{c}z_2\\z_1\end{array}\right).
\end{equation*}
Then the following theorem can be proved: 
\begin{theorem}\label{resavgthm}
Let $z^\eps$ be governed by \eqref{eq:RDS}.
\begin{enumerate}
\item If $\kappa>1$, the process $\{z^{\eps}_t\}_{t\in [0,T]}$  converges as $\eps \to 0$ to the deterministic system
\begin{equation}\label{eq:avgeqz1z2}
\dot{z}=\msB(z).
\end{equation}
\item If $\kappa=1$, the process  $\{z^{\eps}_t\}_{t\in [0,T]}$ converges in law as $\eps \to 0$ to the diffusion given by the SDE
\begin{equation}\label{eq:avgeqz1z2DIFFU}
dz_t=\msB(z_t)dt + \frac{\sigmat}{\sqrt{4\mu}}\left(\begin{array}{c}dW_{1,t}\\ dW_{2,t}\end{array}\right),
\end{equation}
where $W_1, W_2$ are independent standard Wiener processes. 
\end{enumerate}
\end{theorem}
\begin{remark}
We explain why a pair of Wiener processes arise in \eqref{eq:avgeqz1z2DIFFU} even though \eqref{eq:RDS} has only one Wiener process.
Note that the diffusion coefficient of $z^\eps$ in \eqref{eq:RDS} is oscillating fast with frequency $\nu/2\eps$. The diffusion part of the infinitesimal generator for $z^\eps$ can be written as (after expanding $e^{tB/\eps}$):
$$\eps^{2(\kappa-1)}\frac12\frac{\sigma^2}{2\mu}\left(\sin^2\left(\frac{\nu t}{2\eps}\right)\frac{\partial^2}{\partial z_1^2}-\sin\left(\frac{\nu t}{\eps}\right)\frac{\partial^2}{\partial z_1z_2}+\cos^2\left(\frac{\nu t}{2\eps}\right)\frac{\partial^2}{\partial z_2^2}\right).$$
Averaging out the fast oscillations we get
$$\eps^{2(\kappa-1)}\frac12\frac{\sigma^2}{4\mu}\left(\frac{\partial^2}{\partial z_1^2}+\frac{\partial^2}{\partial z_2^2}\right),$$
which is the infinitesimal generator for a pair of independent Brownian motions of strength $\sigma^2/4\mu$.
\end{remark}

In section \ref{subsec:detsyslimitBottom} we show that the state space of the deterministic system \eqref{eq:avgeqz1z2} is the union of domains of attraction of three fixed points. The relation between the fixed points of \eqref{eq:avgeqz1z2} and solutions of \eqref{eq:q1q2sde} is also shown. In section \ref{subsec:transitionBottomDofA} we study how the small noise in \eqref{eq:RDS} facilitates the transition between the domains of attraction.

\subsection{The deterministic system \eqref{eq:avgeqz1z2}.}\label{subsec:detsyslimitBottom}
\begin{figure}\label{fig:z1z2phaseport}
\begin{center}
\begin{overpic}[width=0.5\textwidth]{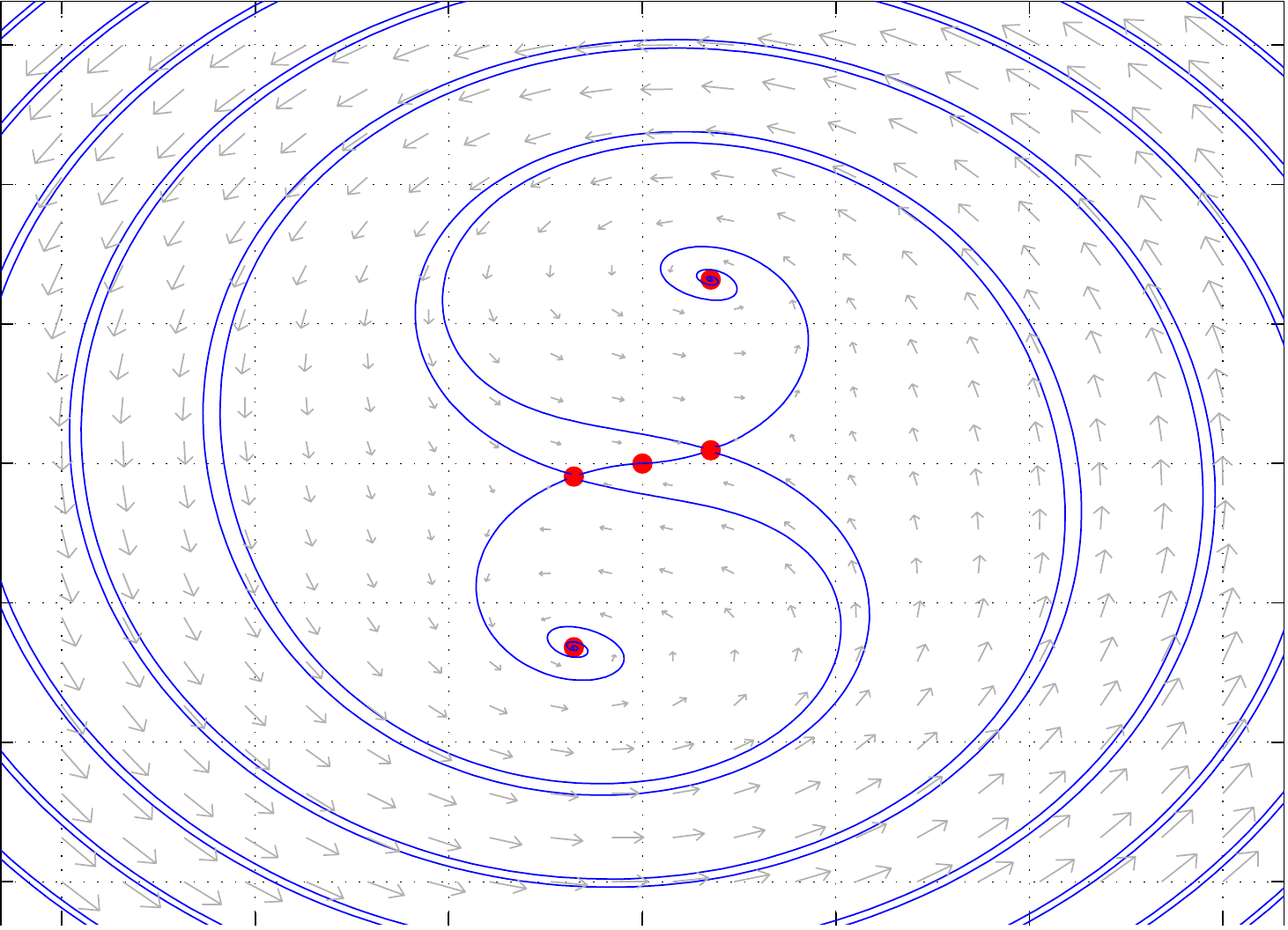}
\put(45,25){$z^{+\pi}$}
\put(53,45){$z^{+0}$}
\put(65,50){$K_1$}
\put(27,24){$K_2$}
\put(40,37){\footnotesize{$K_0$}}
\put(60,37){\footnotesize{$z^{-0}$}}
\end{overpic}
\caption{Typical phase portrait for \eqref{eq:avgeqz1z2} when the conditions \eqref{eq:cond45fixpoints} are satisfied.  The blue lines are stable and unstable manifolds of the saddle points $z^{-0}$ and $z^{-\pi}$. The domain of attraction for $z^0$, $z^{+0}$, $z^{+\pi}$ are separted by the blue lines. Figure generated using the software at \cite{JPolking}.}
\end{center}
\end{figure}
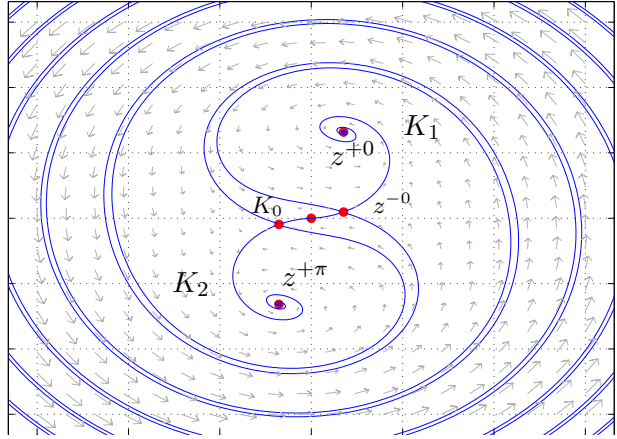

One obvious fixed point of \eqref{eq:avgeqz1z2} is $(0,0)$. Others are given by
\begin{align}\label{eq:saddandperdpoints_a}
\sqrt{z_1^2+z_2^2}&=\frac{\sqrt{4\mu}}{\sqrt{3\gamma}}\sqrt{-(\nu\lambda/2)\pm \sqrt{(\etat\sqrt{2\mu}/4)^2-(\deltat/2)^2}} \quad =:\,{\bf R_{\pm}},\\ \label{eq:saddandperdpoints_b}
\frac{2z_1z_2}{z_1^2+z_2^2}&=\frac{\deltat/2}{\etat \sqrt{2\mu}/4}.
\end{align}
Note that for $\sqrt{z_1^2+z_2^2}$ to be real, we need $\frac14\etat\sqrt{2\mu} > \frac12\deltat.$ 
If 
\begin{equation}\label{eq:cond45fixpoints}
\frac14\etat\sqrt{2\mu} > \frac12\deltat  \qquad \textrm{ and } \qquad  -(\nu\lambda/2)>\sqrt{(\etat\sqrt{2\mu}/4)^2-(\deltat/2)^2},
\end{equation} 
then two values are possible for $\sqrt{z_1^2+z_2^2}$. Also note that if $(z_1,z_2)$ is fixed point then so is $(-z_1,-z_2)$. So, in total there are four nontrivial fixed points. The points with $\sqrt{z_1^2+z_2^2}={\bf R_-}$ are saddles for \eqref{eq:avgeqz1z2} and the points with $\sqrt{z_1^2+z_2^2}={\bf R_+}$ are sinks for \eqref{eq:avgeqz1z2}.
This means that for the system obtained by setting $\sigma=0$ in  \eqref{eq:q1q2sde}, the following solutions are possible (when higher order terms ignored):
\begin{align}\label{eq:origxsolnsposs0}
q_{1,t}^{\eps}&=\sqrt{\mu/\gamma}\, +{\bf{0}}-\eps \frac{\etat\mu \sqrt{\mu/\gamma}}{\nu^2-2\mu}\cos(\nu t), \\ \label{eq:origxsolnspossPa}
q_{1,t}^{\eps}&=\sqrt{\mu/\gamma}\,+ {\bf{\sqrt{\eps}R_+}}\cos(\nu t/2 + \theta^+) -\eps \frac{\etat\mu \sqrt{\mu/\gamma}}{\nu^2-2\mu}\cos(\nu t), \\ \label{eq:origxsolnspossPb}
q_{1,t}^{\eps}&=\sqrt{\mu/\gamma}\,+ {\bf{\sqrt{\eps}R_+}}\cos(\nu t/2 + \theta^+ + \pi) -\eps \frac{\etat\mu \sqrt{\mu/\gamma}}{\nu^2-2\mu}\cos(\nu t), \\ \label{eq:origxsolnspossNa}
q_{1,t}^{\eps}&=\sqrt{\mu/\gamma}\,+ {\bf{\sqrt{\eps}R_-}}\cos(\nu t/2 + \theta^-) -\eps \frac{\etat\mu \sqrt{\mu/\gamma}}{\nu^2-2\mu}\cos(\nu t), \\ \label{eq:origxsolnspossNb}
q_{1,t}^{\eps}&=\sqrt{\mu/\gamma}\,+ {\bf{\sqrt{\eps}R_-}}\cos(\nu t/2 + \theta^- + \pi) -\eps \frac{\etat\mu \sqrt{\mu/\gamma}}{\nu^2-2\mu}\cos(\nu t).
\end{align}
The solutions in \eqref{eq:origxsolnspossNa} and \eqref{eq:origxsolnspossNb} are unstable. The solutions in \eqref{eq:origxsolnsposs0}--\eqref{eq:origxsolnspossPb} are stable. 
Let the fixed points of \eqref{eq:avgeqz1z2} corresponding to \eqref{eq:origxsolnsposs0}--\eqref{eq:origxsolnspossNb} be denoted respectively by $z^0$, $z^{+0}$, $z^{+\pi}$, $z^{-0}$, $z^{-\pi}$. Then, for \eqref{eq:avgeqz1z2}, $z^0$, $z^{+0}, z^{+\pi}$ are stable and $z^{-0}, z^{-\pi}$ are saddles. Let $K_0$ be the domain of attraction  (see figure~\ref{fig:z1z2phaseport}) of the stable trivial equilibrium $(0,0)$. Let $K_1,K_2$ be the domains of attaction of the fixed points $z^{+0}$ and $z^{+\pi}$.
In presence of noise, i.e. $\sigma \neq 0$, transitions occur between the domains of attraction $K_0,\,K_1,\,K_2$ (equivalently between the solutions \eqref{eq:origxsolnsposs0}--\eqref{eq:origxsolnspossPb}).

 Next we study the mean exit time from each of the domains of attraction in the case $\kappa>1$.

\subsection{Transition between the domains of attraction ($\kappa>1$)}\label{subsec:transitionBottomDofA}

We consider the case $\kappa>1$.
Though $z^\eps$ in \eqref{eq:RDS} converges in distribution to the deterministic system \eqref{eq:avgeqz1z2}, the part (ii) of theorem \ref{resavgthm} suggests that, for small $\eps$, the distribution of $z^\eps$ may be better approximated by the distribution of the SDE
\begin{eqnarray}\label{eq:RDSaltapprox}
d\hat{z}_{t}^{\eps} \,=\, \msB(\hat{z}_t^\eps)dt \,+\, \eps^{\kappa-1}\frac{\sigmat}{\sqrt{4\mu}}\left(\begin{array}{c}dW_{1,t}\\ dW_{2,t}\end{array}\right) 
\end{eqnarray}
where $W_1, W_2$ are independent standard Wiener processes. Let $\tau^\eps_i$, for $i=0,1,2$ denote the exit time of \eqref{eq:RDSaltapprox} from the domains $K_i$.  We approximate the mean exit time of $z^\eps$ in \eqref{eq:RDS} by $\expt\, \tau^\eps$.

Mean exit times from a domain $D$ for small noise diffusions of the form $dx=b(x)dt+\eps \sigma dW$ are studied in \cite{FWbook} assuming that the vector field $b$ points inward on the boundary of the domain (see theorem 4.4.1 of \cite{FWbook}). Let $\partial D$ denote the boundary of a domain $D$ and ${\bf n}(z)$ the unit normal vector at point $z$ on $\partial D$. Mean exit times with vector fields $b$ such that $b(z).{\bf n}(z)=0$ are studied in \cite{DayCharBoundary}.

Using theorem 4.2 of \cite{DayCharBoundary} we describe the procedure to get asymptotics of the mean exit time of \eqref{eq:RDSaltapprox} from the domains of attraction $K_i$.

For any $T_1,T_2\in \R$ and $\varphi \in C([T_1,T_2],\R^2)$, define
\begin{equation}\label{E:rate function}
S_{T_1T_2}(\varphi)=\frac12\int_{T_1}^{T_2}\frac{||\dot{\varphi}_t-\msB(\varphi_t)||^2}{\sigma^2/(4\mu)}dt.
\end{equation}
For $i=0,1,2$, define the \emph{quasipotentials}
\begin{equation}\label{eq:quasipotdef}
\msV_i(x):=\inf\{S_{T_1T_2}(\varphi)\,\,:\,\,\varphi\in C([T_1,T_2],\overline{K_i}), \,\,\varphi(T_1)=z_*,\,\,\varphi(T_2)=x,\,\,T_1\leq T_2\},
\end{equation}
where $z_*$ is the fixed point of the domain of attraction $K_i$ with boundary $\partial K_i$. The mean exit time from $K_i$ satisfies
\begin{equation}\label{eq:FWmeanexittimesminV}
\lim_{\eps \to 0}\eps^{2(\kappa-1)}\log \expt \tau_i^\eps \,\,=\,\,\min_{y\in\partial K_i}\msV_i(y).
\end{equation}
The above equations show that the mean exit times $\expt \tau_i^\eps$ are of the order of $e^{{\eps^{-2(\kappa-1)}}}\gg 1$. Though the averaging result that `\eqref{eq:RDSaltapprox} approximates \eqref{eq:RDS}' holds only on times of order $O(1)$, since the fixed points $z_*$ are stable, we expect the mean exit times of \eqref{eq:RDS} to be approximately same exponential order as $\expt \tau_i^\eps$.

Our aim next is to study the dependence of the quasipotentials on the parameters of the system: damping $\delta$, detuning $\lambda$ and strength of nonlinearity $\gamma$ while fixing the values of $\mu$ and strength of periodic excitation $\eta$. The numerical method used to compute the quasipotentials is described in~\ref{sec:app:computingquasipot}.

%############################
\subsection{Dependence of the quasipotentials on the system parameters}\label{subsec:bottomDependQuasipotSysPar}

Recall from \eqref{eq:FWmeanexittimesminV} that the mean exit times from the domain of attraction $K_i$ is determined by $\min_{y\in\partial K_i}\msV_i(y)$. So we need to find $\min_{y\in\partial K_i}\msV_i(y)$. The numerical procedure in~\ref{sec:app:computingquasipot} can be used to find $\msV_i(y)$ and the minmizer can be obtained by inspection. Numerical simulations show that the saddle point on the boundary is the minimizer. This is in agreement with theorem 4.1 in \cite{Day2} which states that saddle points on the boundary are local minimizers for quasipotentials. Recall that $z^{-0}$ and $z^{-\pi}$ are saddles. Define $$V_{0}:=\msV_0(z^{-0}), \quad  V_1:=\msV_1(z^{-0}), \quad  V_2:=\msV_2(z^{-\pi}).$$ Then, mean exit time from $K_i$, is of the order of $e^{V_i/\eps^{2(\kappa-1)}}$. By symmetry $V_1=V_2$. Explicit formulas for $V_0$ and $V_1$ could not be found. However, using numerical simulations some properties of them can be deduced as follows.

 We fix $(\mu,\etat)$ and study how $V_0$ and $V_1$ vary with $(\deltat,\lambda,\gamma)$. We focus only in the regime where there are 5 fixed points for $z$, i.e. the portrait looks as in the figure \ref{fig:z1z2phaseport}. For this situation we need according to \eqref{eq:saddandperdpoints_a} that:
\begin{equation}\label{eq:reqthat5fixpoints}
\deltat\in [0,\sqrt{2\mu}\etat/2], \qquad \lambda < -\frac{1}{\sqrt{2\mu}}\sqrt{(\etat\sqrt{2\mu}/4)^2-(\deltat/2)^2}.
\end{equation} 
Since we fix $(\mu,\etat)$, we can simplify equations by rescaling parameters: $\deltahat=\deltat/(\sqrt{2\mu}\etat/2)$, $\lambdahat=\lambda/(\etat/4)$, $\gammahat=\gamma \frac{3/(4\mu)}{\sqrt{2\mu}\,\etat/4}$. 
Then \eqref{eq:reqthat5fixpoints} becomes $\deltahat \in [0,1]$, and $\lambdahat < -\sqrt{1-\deltah^2}.$
We then get, using $\nu \approx 2\sqrt{2\mu}$, that
\begin{equation}\label{eq:avgdrimsBaltered}
\msB(z)=\frac{\sqrt{2\mu}\etat}{4}\left\{\left(-\gammahat(z_1^2+z_2^2)-\lambdahat \right)\left(\begin{array}{c}z_2\\ -z_1\end{array}\right)-\deltahat\left(\begin{array}{c}z_1\\z_2\end{array}\right)+\left(\begin{array}{c}z_2\\z_1\end{array}\right)\right\}.
\end{equation}
For the fixed points we derive from \eqref{eq:saddandperdpoints_a} and \eqref{eq:saddandperdpoints_b} that
\begin{align}\label{eq:saddandperdpoints_a_norm}
\sqrt{z_1^2+z_2^2}&=\frac{1}{\sqrt{\gammahat}}\sqrt{-\lambdahat \pm \sqrt{1-\deltahat^2}}\qquad =:R_{\pm},\\ \label{eq:saddandperdpoints_b_norm}
\frac{2z_1z_2}{z_1^2+z_2^2}&=\deltahat.
\end{align}
The above two equations suggest that as $(-\lambdahat)$ increases, the size of the domain of attraction of $z^0$ also increases. This can be explained as follows. When $\deltahat$ is fixed, from \eqref{eq:saddandperdpoints_b_norm} we can see that the angle between the `line joining the fixed points to the origin' and the axis $z_1$ is fixed. From \eqref{eq:saddandperdpoints_a_norm} we can see that as $(-\lambdahat)$ increases $R_{\pm}$ increases. A larger $R_{-}$ results in increased size of the domain of attraction of $z^0$.

Intuitively we can deduce one trend: when the damping $\deltahat$ is increased  $V_0$ must increase because a larger damping makes it difficult for the system to reach large values.

Intuition deserts us to predict other dependences. We resort to numerical simulation. We fix $\mu=1$ and $\etat=2$. 
The table below contains the approximate $(V_{0},V_{1})$ pairs obtained by the numerical procedure outlined in \ref{sec:app:computingquasipot}: vertical axis is $\deltahat$ and horizontal axis is $-\lambdahat$.

\begin{table}[h]
    \centering
    \begin{tabular}{|c||c|c|c|c|c|c||}
        \hline
	\backslashbox{$\deltahat$}{$-\lambdahat$} & 1.08 & 1.2 & 1.8 & 2.5 & 5 & 10 \\ \hline 
        0.2 & (0.001, 0.17) & (0.004, 0.17) & (0.04, 0.19) & (0.11, 0.21) & \cellcolor[gray]{0.8}(0.42, 0.23) & \cellcolor[gray]{0.8}(1.18, 0.22)\\\hline
        0.4 & (0.005, 0.29) & (0.014, 0.29) & (0.24, 0.31) & (0.1, 0.29) & \cellcolor[gray]{0.8}(0.88, 0.33) & \cellcolor[gray]{0.8}(2.5, 0.3)\\\hline
	0.6 & (0.022, 0.29) & (0.04, 0.29) & (0.18, 0.28) & \cellcolor[gray]{0.8}(0.4, 0.28) & \cellcolor[gray]{0.8}(1.42, 0.28) & \cellcolor[gray]{0.8}(3.8, 0.3)\\\hline
	0.8 & (0.06, 0.14) & (0.09, 0.14) & \cellcolor[gray]{0.8}(0.33, 0.14) & \cellcolor[gray]{0.8}(0.63, 0.14) & \cellcolor[gray]{0.8}(2.1, 0.14) & \cellcolor[gray]{0.8}(5.2, 0.14)\\\hline
    \end{tabular}
\caption{($V_0,V_1$) values for different values of $\deltahat$ and $\lambdahat$. Other parameters $\mu$ and $\eta$ are fixed at $\mu=1$ and $\etat=2$.}
\end{table}

As expected, as the dissipation $\deltahat$ increases,  $V_{0}$   increases. Whereas, $V_{1}$ increases and then decreases. 

$V_{0}$ increases as $(-\lambdahat)$ increases---possibly because the distance of the saddle from the origin increases. Whereas there is not much variability for $V_{1}$---possibly because $-\gammahat(z_1^2+z_2^2)-\lambdahat \approx -\sqrt{1-\deltahat^2}$ near the $R_+$ fixed point and this is independent of $\lambdahat$; so there is not much $\lambdahat$ dependence in \eqref{eq:avgdrimsBaltered} when $z$ is close to $R_+$.

In the shaded region $V_{0}>V_{1}$, so $q_{t}^{\eps}$ of \eqref{eq:q1q2sde_auxX} spends most of the time close to solution \eqref{eq:origxsolnsposs0}. The unshaded region has $V_{1}>V_{0}$ and so $q_{t}^{\eps}$ of \eqref{eq:q1q2sde_auxX} spends most of the time close to the solutions \eqref{eq:origxsolnspossPa}-\eqref{eq:origxsolnspossPb}.

\subsection*{Summary} 
The capture of an oscillatory nonlinear system into resonance by periodic perturbations is an important process in many applications. When noise perturbations are also present, the noise facilitates escape from resonance zone.
This paper used averaging techniques to determine the effects of noisy perturbations on the passage of trajectories through a resonance zone.
We have examined a prototypical beam type nonlinear energy harvesting model that contains a double well potential presented in ~\cite{McInnes2008,Cartmell2014} and have shown that the averaging technique enables some of the basic behavior of the model to be simply determined.  In particular, we have obtained the mean exit time for the pendulum Hamiltonian $\mathcal{H}^\eps$ which governs the rate of escape from the resonance zone.

\subsection*{Acknowledgements}
The authors would like to acknowledge the support of the AFOSR under grant number FA9550-12-1-0390  and the National Science Foundation~(NSF) under grant number CMMI 1030144. 
Any opinions, findings, and conclusions or recommendations expressed in this paper are those of the authors and do not necessarily reflect the views of the NSF.

%#####################################
%################ APPENDIX ###########
%#####################################

\appendix

\section{Computation of $I_r$ and $J_r$}\label{appsec:computeIJ}
The system \eqref{eq:q1q2sde} with $\eps=0$ is the same as the system defined by \eqref{eq:unpertsysHamsys}-\eqref{eq:PotDefprimary}. We fix $\mu=1$, $\gamma=1$. Then $H$ can take values in $[-1/4,\infty)$. The values $[-1/4,0)$ correspond to the region inside the homoclinic orbits (see figure \ref{fig:potandcontourplot}), $H=0$ denotes the homoclinic orbit, and the region outside the homoclinic orbit corresponds to $H>0$.

Solution to \eqref{eq:unpertsysHamsys}-\eqref{eq:PotDefprimary} for each level of the Hamiltonian can be written in terms of elliptic functions \cite{ByrdFriedman}. 
So we use elliptic modulus $k$ as a proxy for the Hamiltonian $H$.

In the following, $\jsn, \jcn, \jdn$ are Jacobi elliptic functions, and $K, E$ are complete elliptic integrals of the first and second kind respectively. In all the elliptic functions, the elliptic modulus used is $k$.

\subsection{The case of $-1/4\leq H <0$}
Given $H$, let $k$ be defined by 
$$H=-\frac{1-k^2}{(2-k^2)^2}.$$
The solution to \eqref{eq:unpertsysHamsys}-\eqref{eq:PotDefprimary} can be written as
$$q_1(\varphi)=\sqrt{\frac{2}{2-k^2}}\jdn(\frac{2K}{2\pi}\varphi), \qquad q_2(\varphi)=-\frac{\sqrt{2}k^2}{2-k^2}\jcn(\frac{2K}{2\pi}\varphi)\jsn(\frac{2K}{2\pi}\varphi), \qquad \dot{\varphi}=\Omg,$$
where $\Omg=\frac{2\pi}{2K\sqrt{2-k^2}}$. The action $I$ by definition is $\frac{1}{2\pi}\int q_2\,dq_1$ with the integral over the Hamiltonian orbit, and $I$ turns out to be
$$I=\frac{2}{3\pi(2-k^2)^{3/2}}\bigg(2(k^2-1)K+(2-k^2)E\bigg).$$
Now, $m:n$ resonance occurs when $k$ is such that we have $\nu = \Omg\frac{m}{n}$. Let $I_r$ denote the value of $I$ at this resonance.
Noting that $\frpt{I}{q_2}=\frpt{I}{H}\frpt{H}{q_2}=\frac{1}{\Omg}q_2$ we get that
$$\mfF(I,\varphi,\theta)=\frac{q_2(I,\varphi)}{\Omg(I)}(q_1(I,\varphi)\eta\cos\theta + \alpha \cos\theta -\delta q_2(I,\varphi)).$$
Evaluating $\la\mfF(\psi)\ra=\frac{1}{2m\pi}\int_0^{2m\pi}\mfF(I,\psi+\frac{n}{m}\theta,\theta)d\theta$ with $I=I_r$ gives $$\la \mfF(\psi)\ra = -\delta I_r+J_r \sin(m\psi/n)= -\delta I_r+ (\eta J^\eta_r + \alpha J^\alpha_r)\sin(m\psi/n)$$
with  
$$J^\eta_r=-\frac{\pi^2(m/n)^2}{2K^2(2-k^2)}{\tt csch}\left(\frac{m}{n}\frac{\pi \widetilde{K}}{K}\right)\mathbf{1}_{\{m/n\in \mathbb{N}\}},$$
$$J^\alpha_r=-\frac{\pi (m/n)}{\sqrt{2}K\sqrt{2-k^2}}{\tt sech}\left(\frac{m}{n}\frac{\pi \widetilde{K}}{K}\right)\mathbf{1}_{\{m/n\in \mathbb{N}\}},$$
where $\widetilde{K}$ is the complete elliptic integral of the first kind with modulus $\sqrt{1-k^2}$.
To calculate $\la \mfG(\psi)\ra$ we make note of the following relation (see lemma 3.4 of \cite{YagasakiSIAM96}):
$$ \frac{\partial}{\partial \varphi}\mfG(I,\varphi,\theta)+ \frac{\partial}{\partial I}\mfF(I,\varphi,\theta)=-\delta.$$
Hence
$$ \frac{\partial}{\partial \psi}\mfG(I,\psi+\frac{n}{m}\theta,\theta)+ \frac{\partial}{\partial I}\mfF(I,\psi+\frac{n}{m}\theta,\theta)=-\delta.$$
Now,
\begin{align*}
\frac{\partial}{\partial \psi}\la \mfG(\psi)\ra &= \frac{\partial}{\partial \psi}\frac{1}{2m\pi}\int_0^{2m\pi}\mfG(I,\psi+\frac{n}{m}\theta,\theta)d\theta \\
&=\frac{1}{2m\pi}\int_0^{2m\pi}\frac{\partial}{\partial \psi}\mfG(I,\psi+\frac{n}{m}\theta,\theta)d\theta \\
&=-\delta -\frac{\partial}{\partial I}\frac{1}{2m\pi}\int_0^{2m\pi}\mfF(I,\psi+\frac{n}{m}\theta,\theta)d\theta \\
&=-\delta - \left(-\delta + (\eta (J^{m:n}_{\eta})'+\alpha (J^{m:n}_{\alpha})')\sin(m\psi /n)\right) \\
&=-(\eta (J^{\eta}_r)'+\alpha (J^{\alpha}_r)')\sin(m\psi /n)
\end{align*}
which gives on integrating 
$$\la \mfG(\psi)\ra = \frac{n}{m}(\eta (J^{\eta}_r)'+\alpha (J^{\alpha}_r)')\cos(m\psi /n),$$ 
where $'$ indicates differentiation w.r.t $I$. These derivatives can be evaluated as $(J^{\eta}_r)'=\frac{\partial J^{\eta}_r}{\partial k}/\frac{\partial I}{\partial k}$ etc.

\subsection{The case of $H>0$}
Given $H$, let $k$ be defined by 
$$H=\frac{k^2(1-k^2)}{(2k^2-1)^2}.$$
The solution to \eqref{eq:unpertsysHamsys}-\eqref{eq:PotDefprimary} can be written as
$$q_1(\varphi)=\sqrt{\frac{2k^2}{2k^2-1}}\jcn(\frac{4K}{2\pi}\varphi), \qquad q_2(\varphi)=-\frac{\sqrt{2}k}{2k^2-1}\jsn(\frac{4K}{2\pi}\varphi)\jdn(\frac{4K}{2\pi}\varphi), \qquad \dot{\varphi}=\Omg,$$
where $\Omg=\frac{2\pi}{4K\sqrt{2k^2-1}}$. The action $I$ turns out to be
$$I=\frac{4}{3\pi(2k^2-1)^{3/2}}\bigg((1-k^2)K+(2k^2-1)E\bigg).$$
Using similar procedure as in the case of $H\in [-1/4,0)$, we get
$$J^\eta_r=-\frac{\pi^2(m/n)^2}{4K^2(2k^2-1)}{\tt csch}\left(\frac{m}{2n}\frac{\pi \widetilde{K}}{K}\right)\mathbf{1}_{\{\frac{m}{n}\in 2\mathbb{N}\}},$$
$$J^\alpha_r=-\frac{\pi (m/n)}{\sqrt{2}K\sqrt{2k^2-1}}{\tt sech}\left(\frac{m}{2n}\frac{\pi \widetilde{K}}{K}\right)\mathbf{1}_{\{\frac{m}{n}\in 2\mathbb{N}-1\}},$$
$$\la \mfG(\psi)\ra = \frac{n}{m}(\eta (J^{\eta}_r)'+\alpha (J^{\alpha}_r)')\cos(m\psi /n).$$ 

%##############################
%##############################appsec:showingB1eq0
%##############################
\section{Showing $\mfB_1 \equiv 0$.}\label{appsec:showingB1eq0}
Let $$\bbK_0(\varphi)=\frac{1}{\Omg}\bigg(\eta q_1(I_r,\varphi)q_2(I_r,\varphi)+\alpha q_2(I_r,\varphi)\bigg),$$
$$\bbK_c(\varphi)=\bbK_0(\varphi)\cos(m\varphi/n), \qquad \bbK_s(\varphi)=\bbK_0(\varphi)\sin(m\varphi/n), \qquad \bbK_\delta(\varphi)=-\frac{1}{\Omg}\delta q_2^2(I_r,\varphi),$$
with $q_i(I,\varphi)$ being represented in terms of elliptic functions with modulus $k$ and argument $\frac{4K(k)}{2\pi}\varphi$ where $K(k)$ is the complete elliptic integral of the first kind.
Then,
$$\mfF(\psih+\frac{n}{m}\theta,\theta)=\cos(\frac{m\psih}{n})\bbK_c(\psih+\frac{n}{m}\theta)+\sin(\frac{m\psih}{n})\bbK_s(\psih+\frac{n}{m}\theta)+\bbK_\delta(\psih+\frac{n}{m}\theta).$$
It can be verified easily that
$$\nu \mathcal{A}_1(\psih,\theta)\,\,=\,\,-\la \mfF(\psih)\ra\theta\,\, +\,\, \frac{m}{n}\int_{\psih}^{\psih+\frac{n}{m}\theta}\left(\cos(\frac{m\psih}{n})\bbK_c(\varphi)+\sin(\frac{m\psih}{n})\bbK_s(\varphi)+\bbK_\delta(\varphi)\right)d\varphi.$$

The average w.r.t $\theta$ of the function $\theta \mapsto \mathcal{A}_1(\psih,\theta)\mfF(\psih+\frac{n}{m}\theta,\theta)$ is given by
\begin{align*}
\la \mathcal{A}_1\mfF\ra &= \frac{1}{2n\pi \nu}\int_{\psih}^{\psih+2n\pi}\left\{-\la \mfF(\psih)\ra\frac{(\hat{\varphi}-\psih)m}{n} +\right. \\
 &\qquad \left. + \frac{m}{n}\int_{\psih}^{\hat{\varphi}}\left(\cos(\frac{m\psih}{n})\bbK_c(\varphi)+\sin(\frac{m\psih}{n})\bbK_s(\varphi)+\bbK_\delta(\varphi)\right)d\varphi\right\}\times \\
 &\qquad \qquad \times \left\{\cos(\frac{m\psih}{n})\bbK_c(\hat{\varphi})+\sin(\frac{m\psih}{n})\bbK_s(\hat{\varphi})+\bbK_\delta(\hat{\varphi}) \right\}d\hat{\varphi}.
\end{align*}
\begin{align*}
\nu \frpt{\mathcal{A}_1}{\psih}(\psih,\theta)\,\,=\,\,-\frpt{\la \mfF(\psih)\ra}{\psih}\theta\,\, &+\,\, \frac{m}{n}\int_{\psih}^{\psih+\frac{n}{m}\theta}\left(-\frac{m}{n}\sin(\frac{m\psih}{n})\bbK_c(\varphi)+\frac{m}{n}\cos(\frac{m\psih}{n})\bbK_s(\varphi)\right)d\varphi \\
 \qquad \qquad &+ \frac{m}{n}\cos(\frac{m\psih}{n})(\bbK_c(\psih+\frac{n}{m}\theta)-\bbK_c(\psih)) \\
 \qquad \qquad &+ \frac{m}{n}\sin(\frac{m\psih}{n})(\bbK_s(\psih+\frac{n}{m}\theta)-\bbK_s(\psih)) \\
 \qquad \qquad &+ \frac{m}{n}(\bbK_\delta(\psih+\frac{n}{m}\theta)-\bbK_\delta(\psih)).
\end{align*}
Hence  
\begin{align*}
\la \frpt{\mathcal{A}_1}{\psih}\ra = \frac{1}{2n\pi \nu}\int_{\psih}^{\psih+2n\pi}\bigg\{\frac{m}{n}\int_{\psih}^{\hat{\varphi}}&\left(-\frac{m}{n}\sin(\frac{m\psih}{n})\bbK_c(\varphi)+\frac{m}{n}\cos(\frac{m\psih}{n})\bbK_s(\varphi)\right)d\varphi \,+\,\\
 \qquad \qquad & +\,\frac{m}{n}\cos(\frac{m\psih}{n})(\bbK_c(\hat{\varphi})-\bbK_c(\psih))  \\
 \qquad \qquad & +\,\frac{m}{n}\sin(\frac{m\psih}{n})(\bbK_s(\hat{\varphi})-\bbK_s(\psih))  \\
 \qquad \qquad & +\,\frac{m}{n}(\bbK_\delta(\hat{\varphi})-\bbK_\delta(\psih)) \,\,-\,\,\frac{m}{n}\frpt{\la \mfF(\psih)\ra}{\psih}(\hat{\varphi}-\psih)\bigg\} d\hat{\varphi}.
\end{align*}

We start with the above and using integration-by-parts show that $\mfB_1\equiv 0$. It is convenient to introduce some notation first:
For periodic functions with period $2n\pi$ we define
\begin{align}
\{f\}&=\frac{1}{2n\pi}\int_{\psi}^{\psi+2n\pi}f(\varphi)d\varphi, \\
\{f,g\}&=\frac{1}{2n\pi}\int_{\psi}^{\psi+2n\pi}g(\hat{\varphi})\left(\int_{\psi}^{\hat{\varphi}}f(\varphi)d\varphi \right)d\hat{\varphi},\\
\{\{f\}\}&=\{1,f\}.
\end{align}
Note that $\{f\}$ does not depend on $\psi$, but $\{f,g\}$ and $\{\{f\}\}$ do. Further, 
\begin{equation}\label{eq:auxpropfggf}
\{f,g\}+\{g,f\}=2n\pi\{f\}\{g\},
\end{equation}
and
\begin{equation}\label{eq:auxpropfggf2}
\frac{1}{2n\pi}\int_{\psi}^{\psi+2n\pi}(\varphi-\psi)f(\varphi)d\varphi = 2n\pi\{f\}-\{\{f\}\}.
\end{equation}

Let $\mfc=\cos(m\psi/n)$ and $\mfs=\sin(m\psi/n)$. 

It can be verified that  
$$\{\bbK_\delta\}=-\delta I_r, \qquad \{\bbK_s\}=J_r, \qquad \{\bbK_c\}=0.$$
Hence
\begin{equation}\label{eq:auxpropfggf3}
\la \mfF \ra = \{\bbK_\delta\}+\{\bbK_s\}\mfs.
\end{equation}

Akin to the results \eqref{eq:intbypartAux1} and \eqref{eq:intbypartAux2}, we can prove the following lemma:
\begin{lemma}\label{lem:Aavgauxresults_2} The following five identities hold:
\begin{enumerate}
\item $\Aavg[\{\{\bbK_c\}\}\mfc(\{\bbK_\delta\}+\{\bbK_s\}\mfs)]=\Aavg[\{\{\bbK_c\}\}\frac{m}{n}\Omg'_rh^2\mfs]-\Aavg[\{\{\bbK_c\}\}'\Omg'_rh^2\mfc]$, 
\item $\{\bbK_s\}\Aavg[\{\{\bbK_s\}\}\mfs^2]= -\{\bbK_\delta\}\Aavg[\{\{\bbK_s\}\}\mfs]-\Aavg[\{\{\bbK_s\}\}'\Omg'_rh^2\mfs]-\Aavg[\{\{\bbK_s\}\}\frac{m}{n}\Omg'_rh^2\mfc]$,
\item $-\frac{2n\pi}{2}\{\bbK_s\}^2\Aavg[\mfs^2]=-\frac{2n\pi}{2}\{\bbK_\delta\}^2+\frac{2n\pi}{2}\{\bbK_s\}\frac{m}{n}\Omg'_r\Aavg[h^2\mfc]$, 
\item $\{\bbK_\delta\}\Aavg[\{\{\bbK_s\}\}] + \Aavg[\Omg'_rh^2\{\{\bbK_s\}\}'] +\{\bbK_s\}\Aavg[\{\{\bbK_s\}\}\mfs] =0$,
\end{enumerate}
and $(v)$
\begin{align*}\Aavg[\mfs( \{\bbK_\delta\}\{\{\bbK_s\}\}+\{\bbK_s\}\{\{\bbK_\delta\}\}-2n\pi\{\bbK_\delta\}\{\bbK_s\})]= -\frac{\{\bbK_\delta\}^2}{\{\bbK_s\}}\Aavg[\{\{\bbK_s\}\}]-\{\bbK_\delta\}\Aavg[\{\{\bbK_\delta\}\}]  \\
 \qquad +2n\pi\{\bbK_\delta\}^2 -\frac{\{\bbK_\delta\}}{\{\bbK_s\}}\Aavg[\Omg'_rh^2\{\{\bbK_s\}\}']  -\Aavg[\Omg'_rh^2\{\{\bbK_\delta\}\}'].
\end{align*}
\end{lemma}

Now,
\begin{align*}
\frac{n}{m}\nu \la \mathcal{A}_1\mfF\ra &= \mfc(\{\bbK_\delta,\bbK_c\}+\{\bbK_c,\bbK_\delta\}) + \mfs(\{\bbK_\delta,\bbK_s\}+\{\bbK_s,\bbK_\delta\}) \\
 \qquad \qquad & \qquad  + \mfc^2\{\bbK_c,\bbK_c\} + \mfs^2\{\bbK_s,\bbK_s\} + \mfc\mfs(\{\bbK_c,\bbK_s\}+\{\bbK_s,\bbK_c\}) \\
 \qquad \qquad & \qquad + \{\bbK_\delta,\bbK_\delta\} - \la \mfF\ra\mfc (2n\pi\{\bbK_c\}-\{\{\bbK_c\}\}) \\
 \qquad \qquad & \qquad - \la \mfF\ra\mfs (2n\pi\{\bbK_s\}-\{\{\bbK_s\}\}) - \la \mfF\ra (2n\pi\{\bbK_\delta\}-\{\{\bbK_\delta\}\}).
\end{align*}
Using \eqref{eq:auxpropfggf}, \eqref{eq:auxpropfggf3} and that $\{\bbK_c\}=0$ we have
\begin{align*}
\frac{n}{m}\nu \la \mathcal{A}_1\mfF\ra &= -\frac{2n\pi}{2}\{\bbK_\delta\}^2+\{\bbK_\delta\}\{\{\bbK_\delta\}\} + \{\{\bbK_c\}\}\mfc(\{\bbK_\delta\}+\{\bbK_s\}\mfs) \\
 \qquad \qquad &\qquad + \mfs ( \{\bbK_\delta\}\{\{\bbK_s\}\}+\{\bbK_s\}\{\{\bbK_\delta\}\}-2n\pi\{\bbK_\delta\}\{\bbK_s\})\\
\qquad \qquad & \qquad + \mfs^2(-\frac{2n\pi}{2}\{\bbK_s\}^2+\{\bbK_s\}\{\{\bbK_s\}\}).
\end{align*}
Using lemma \ref{lem:Aavgauxresults_2} we have
\begin{align*}
\frac{n}{m}\nu \Aavg[\la \mathcal{A}_1\mfF\ra] &= -\frac{2n\pi}{2}\{\bbK_\delta\}^2+\{\bbK_\delta\}\Aavg\{\{\bbK_\delta\}\}  + \Aavg[\{\{\bbK_c\}\}\frac{m}{n}\Omg'_rh^2\mfs]\\
\qquad \qquad & \qquad -\Aavg[\{\{\bbK_c\}\}'\Omg'_rh^2\mfc] -\frac{\{\bbK_\delta\}^2}{\{\bbK_s\}}\Aavg[\{\{\bbK_s\}\}]-\{\bbK_\delta\}\Aavg[\{\{\bbK_\delta\}\}]\\
 \qquad \qquad & \qquad +2n\pi\{\bbK_\delta\}^2  -\frac{\{\bbK_\delta\}}{\{\bbK_s\}}\Aavg[\Omg'_rh^2\{\{\bbK_s\}\}']  -\Aavg[\Omg'_rh^2\{\{\bbK_\delta\}\}']\\
\qquad \qquad & \qquad -\frac{2n\pi}{2}\{\bbK_\delta\}^2+\frac{2n\pi}{2}\{\bbK_s\}\frac{m}{n}\Omg'_r\Aavg[h^2\mfc] -\{\bbK_\delta\}\Aavg[\{\{\bbK_s\}\}\mfs] \\
 \qquad \qquad& \qquad -\Aavg[\{\{\bbK_s\}\}'\Omg'_rh^2\mfs]-\Aavg[\{\{\bbK_s\}\}\frac{m}{n}\Omg'_rh^2\mfc].
\end{align*}
Simplifying and rearranging we have
\begin{align*}
\frac{n}{m}\nu \Aavg[\la \mathcal{A}_1\mfF\ra] &=  \Aavg[(\mfs\{\{\bbK_c\}\}-\mfc\{\{\bbK_s\}\})\frac{m}{n}\Omg'_rh^2]+\frac{2n\pi}{2}\{\bbK_s\}\frac{m}{n}\Omg'_r\Aavg[h^2\mfc]\\
 &\qquad  -\Aavg[(\{\{\bbK_c\}\}'\mfc+\{\{\bbK_s\}\}'\mfs + \{\{\bbK_\delta\}\}') \Omg'_rh^2] \\
 &\qquad  -\frac{\{\bbK_\delta\}^2}{\{\bbK_s\}}\Aavg[\{\{\bbK_s\}\}]  -\frac{\{\bbK_\delta\}}{\{\bbK_s\}}\Aavg[\Omg'_rh^2\{\{\bbK_s\}\}'] -\{\bbK_\delta\}\Aavg[\{\{\bbK_s\}\}\mfs]. 
\end{align*}
Using that $\{\{\bbK_s\}\}'=-\bbK_s(\psi)+\{\bbK_s\}$ etc we have
$$
\{\{\bbK_c\}\}'\mfc+\{\{\bbK_s\}\}'\mfs + \{\{\bbK_\delta\}\}' = (\{\bbK_c\}\mfc + \{\bbK_s\}\mfs + \{\bbK_\delta\})  - (\bbK_c(\psi)\mfc + \bbK_s(\psi)\mfs + \bbK_\delta(\psi)),
$$
and employing this above we have
\begin{align*}
\frac{n}{m}\nu \Aavg[\la \mathcal{A}_1\mfF\ra] =  \Aavg[(\mfs\{\{\bbK_c\}\}-\mfc\{\{\bbK_s\}\})\frac{m}{n}\Omg'_rh^2]  -\Aavg[(\{\bbK_c\}\mfc + \{\bbK_s\}\mfs + \{\bbK_\delta\}) \Omg'_rh^2] \\
 \qquad \quad +\Aavg[(\bbK_c(\psi)\mfc + \bbK_s(\psi)\mfs + \bbK_\delta(\psi)) \Omg'_rh^2] +\frac{2n\pi}{2}\{\bbK_s\}\frac{m}{n}\Omg'_r\Aavg[h^2\mfc]. 
\end{align*}
Inspection of $\frac{n}{m}\nu \Aavg[\Omg'h^2\la \frpt{\mathcal{A}_1}{\psi}\ra]$ shows that exactly the above terms arise but with opposite sign.

%##################
%##################
%##################

%##################
%##################
%##################

%###################
%###################
%###################

\section{Computation of the quasipotential}\label{sec:app:computingquasipot}
Recall the definition \eqref{eq:quasipotdef} of the quasipotential $\msV$ and the action functional $S_{T_1T_2}$ defined in \eqref{E:rate function}.
The optimization problem in  \eqref{eq:quasipotdef} can be written as follows:
$$
\textrm{Minimize }  \frac{1}{2}\int_{T_1}^{T_2}||u_s||_2^2ds  \quad \textrm{ subject to } \quad \dot{z}_t=\mfB(z_t)+u_t  \quad \textrm{ with } \quad z_{T_1}=z_*, \quad z_{T_2}=x.
$$
Note that $T_1$ and $T_2$ are also free in the optimization, i.e. the minimum is over all possible $T_1, T_2$ with $T_1\leq T_2$.

The usual method to solve this optimal control problem is as follows: 

Define the Hamiltonian
\begin{equation}\label{eq:optcontHamdef}
H(z,p):= \sup_{u} \left(p^{tr}(\mfB(z)+u)-\frac12||u||_2^2\right).
\end{equation}
It is easy to see that the sup is obtained by taking $u=p$ and so 
\begin{equation}\label{eq:optcontHamdefeval}
H(z,p)=p^{tr}\mfB(z)+ \frac12||p||_2^2.
\end{equation}
Then the trajectories for which $\frac{1}{2}\int_{T_1}^{T_2}||u_s||_2^2ds$ has first variation zero satisfy the Euler-Lagrange equations
\begin{equation}\label{eq:optcontHameqEL}
\dot{z}=\frac{\partial H}{\partial p}, \qquad \dot{p}=-\frac{\partial H}{\partial z}.
\end{equation}
Further, the fact that the time variables $T_i$ are free, forces 
\begin{equation}\label{eq:optcontFreeTimpHequiv0}
H\equiv 0.
\end{equation}
If we are interested in calculating the quasipotential $\msV_i(x)$, we need to impose the boundary conditions
\begin{equation}\label{eq:optcont2pointbc}
z(T_1)=z_*, \qquad z(T_2)=x
\end{equation}
where $z_*$ is the stable fixed point whose domain of attraction is $K_i$. Note that these are four boundary conditions---two for $T_1$ and two for $T_2$. The function $\msV_i(x)$ itself is obtained by integrating
\begin{equation}\label{eq:optcontIntg2getquasipot}
\dot{V}=\frac{1}{2(\sigma^2/4\mu)}||p||^2\quad \textrm{ for } t\in [T_1,T_2], \quad V(T_1)=0,
\end{equation}
and setting $\msV_i(x)=V(T_2)$. Note that the sup in \eqref{eq:optcontHamdef} is attained at $u=p$.

Hence, the quasipotential can be obtained by solving \eqref{eq:optcontHameqEL} for the 4-dimensional system $(z,p)$ with the four boundary conditions \eqref{eq:optcont2pointbc} while using \eqref{eq:optcontFreeTimpHequiv0} to determine the free parameters $T_i$ and then using \eqref{eq:optcontIntg2getquasipot}.

The above suggested method works except for the following issue. Recall that the $z_*$ in the definition of quasipotential is a fixed point for $\dot{z}=\msB(z)$. Hence $\msB(z_*)=0$. When $z=z_*$,  \eqref{eq:optcontFreeTimpHequiv0} implies that $p=0$. So $(z,p)=(z_*,0)$ is a fixed point for the system of equations \eqref{eq:optcontHameqEL} and the system started at $(z_*,p)$ does not move from it.

To rectify this, \cite{Day1} suggests the following as a numerical procedure to calculate the quasipotential. The optimization above does not occur for finite times $T_1, T_2$. Optimal trajectory takes infinite time to leave from $(z_*,p)$. When it leaves, it leaves along the unstable manifold of \eqref{eq:optcontHameqEL} at $(z_*,p=0)$. So, instead of starting at $(z_*,0)$, start \eqref{eq:optcontHameqEL} at a point on the unstable manifold but very close to $(z_*,0)$. The unstable manifold at $(z_*,0)$ is tangential to the unstable eigenspace of the linearization of system \eqref{eq:optcontHameqEL} at $(z_*,0)$. And this tangent can be easily found.  Given $z^\dag$ very close to $z_*$, there is a unique $p^\dag$ so that $(z^\dag,p^\dag)$ belongs to the unstable eigenspace. So, we pick lot of $z^\dag$ close to $z_*$ and find corresponding $p^\dag$ and simulate \eqref{eq:optcontHameqEL}. Of all these simulations whichever trajectory passes through $x$ is the desired trajectory.

The above is the numerical procedure that we use to study the dependence of the quasipotential on the system parameters in section \ref{subsec:bottomDependQuasipotSysPar}. For the sake of completeness we write the system \eqref{eq:optcontHameqEL} explicitly, clearly showing its linear and nonlinear parts. Let $c=-\frac{3\gamma}{4\mu}||z_*||_2^2-\frac12\nu\lambda$. Then \eqref{eq:optcontHameqEL} can be written as
\begin{equation}\label{eq:optcontHameqEL_explicit}
\left(\begin{array}{c}\dot{z}_1\\ \dot{z}_2\\ \dot{p}_1\\ \dot{p}_2\end{array}\right)=\left(\begin{array}{cc} M & I_{2\times 2} \\ 0_{2\times 2} & N \end{array}\right)\left(\begin{array}{c}z_1\\z_2\\p_1\\p_2\end{array}\right)+\frac{3\gamma}{4\mu}\left(\begin{array}{c}-(||z||_2^2-||z_*||_2^2)z_2 \\ (||z||_2^2-||z_*||_2^2)z_1 \\  2z_1(p_1z_2-p_2z_1) \\ 2z_2(p_1z_2-p_2z_1)\end{array}\right)
\end{equation}
where $M=\left(\begin{array}{cc}-\delta/2 & c+\eta\sqrt{2\mu}/4  \\ -c+\eta\sqrt{2\mu}/4 & -\delta/2  \end{array}\right)$, and $N=\left(\begin{array}{cc}\delta/2 & c-\eta\sqrt{2\mu}/4  \\ -c-\eta\sqrt{2\mu}/4 & \delta/2  \end{array}\right)$. 
Let $U$ be a matrix such that $(z=Up, p)$ is in the unstable eigenspace of the linearized system. Then we have $I+MU-UN=0$. After solving this for $U$ we can take a point in the unstable eigenspace as $(z,U^{-1}z)$. So, we choose $z$ values near by $z_*$ and then start \eqref{eq:optcontHameqEL_explicit} at $(z,U^{-1}z)$.

%##################
%##################
%##################
\vspace{10pt}
\noindent {\bf \large{References}}


\vspace{10pt}

\begin{thebibliography}{99}
\bibitem{DaqMasSurv}
Daqaq MF, Masana R, Erturk A, Quinn D. 2014. \newblock On the role of nonlinearities in vibratory energy harvesting: A critical review and discussion. {\it ASME. Appl. Mech. Rev.} {\bf 66}(4):040801. \url{http://dx.doi.org/10.1115/1.4026278}.

\bibitem{McInnes2008}
McInnes CR, Gorman DG, Cartmell MP. 2008. \newblock Enhanced vibrational energy harvesting using nonlinear stochastic resonance.
\newblock  \emph{Journal of Sound and Vibration.} {\bf 318}:655-662 \url{http://dx.doi.org/10.1016/j.jsv.2008.07.017}

\bibitem{Cartmell2014}
Zheng R, Nakano K, Hu H, Su D, Cartmell MP. 2014. \newblock An application of stochastic resonance for energy harvesting in a bistable vibrating system.
\newblock \emph{Journal of Sound and Vibration}. {\bf 333}(12):2568-2587. \url{http://dx.doi.org/10.1016/j.jsv.2014.01.020}

\bibitem{Bol64}
{Bolotin VV.} 1964.
\newblock \emph{The dynamic stability of elastic systems}.
\newblock Holden-Day, San Francisco.

\bibitem{Nay79}
{Nayfeh AH, Mook DT.} 1979.
\newblock \emph{Nonlinear oscillations}.
\newblock John Wiley \& Sons, New York.

\bibitem{Guc83}
{Guckenheimer J, Holmes P.} 1983.
\newblock {\em Nonlinear oscillations, dynamical systems, and bifurcations of vector fields}.
\newblock Springer-Verlag, New York.

\bibitem{Morozov1998}
{Morozov AD.} 1998.
\newblock {\em Quasi-conservative systems: cycles, resonances and chaos}.
\newblock World Scientific Series on Nonlinear Science Series A: Volume 30.

\bibitem{Bol84}
{Bolotin  VV.} 1984.
\newblock \emph{Random vibrations of elastic systems}.
\newblock Martinus Nijhoff Publishers, Boston.


\bibitem{Sri:Ebe86}
Ebeling W, Herzel H, Richert W, Schimansky-Geier L. 1986.
\newblock Influence of noise on Duffing-van der Pol oscillations.
\newblock {\em ZAMM}, {\bf 66}(3):141-- 146.

\bibitem{ArnL96}
Arnold L, Namachchivaya N, Schenk KL. 1996.
\newblock Toward an understanding of stochastic {H}opf bifurcations: 
A case study.
\newblock {\em Journal of Bifurcation and Chaos}, {\bf 6}(11):1947--1975.
\url{http://dx.doi.org/10.1142/S0218127496001272 }

\bibitem{Lia99}
Liang Y, Sri Namachchivaya N. 1999.
\newblock P-bifurcation in the stochastic version of the {D}uffing-van der {P}ol equation.
\newblock In H.~Crauel and V.~M. Gundlach, editors, {\em Stochastic Dynamics},
  pages 146--170. Springer-Verlag, Berlin.

\bibitem{SriSowVed}
Sri Namachchivaya N, Sowers RB, Vedula L. 2001.
\newblock Nonstandard reduction of noisy Duffing-van der Pol equation.
\newblock \emph{Dynamical Systems}. {\bf 16}(3):201-224.
  \url{http://dx.doi.org/10.1080/14689360118168}
  
\bibitem{Sri:Fre98}
Freidlin MI, Weber M. 1998.
\newblock Random perturbations of nonlinear oscillators.
\newblock {\em The Annals of Probability}, 26(3):925--967. \url{https://projecteuclid.org/euclid.aop/1022855739}

\bibitem{Sri2001b}
Sri Namachchivaya N, Sowers RB. 2001.
\newblock Unified approach for noisy nonlinear Mathieu-type systems.
\newblock \emph{Stochastics \& Dynamics}. {\bf 1}(3):405-450.
\url{http://dx.doi.org/10.1142/S0219493701000217 }

\bibitem{ArnoldDynsys}
{Arnold VI. (Ed.)} 1987. {\em Dynamical systems III}. Springer-Verlag.

\bibitem{FWMult}
{Freidlin MI, Wentzell AD.} 2003. \newblock{Averaging principle for stochastic perturbations of multifrequency systems.} \emph{Stochastics and Dynamics}, {\bf 3}:393-408. \url{http://dx.doi.org/10.1142/S0219493703000747}

\bibitem{Neishtadt1991}
Neishtadt AI. 1991.
\newblock Averaging and passage through resonances.
\newblock \emph{Proceedings International Congress Mathematicians}, 
The Mathematical Society of Japan, 1271-1284.

\bibitem{SowersLDP}
{Freidlin MI, Sowers RB.} 1999. \newblock{A comparison of homogenization and large deviations, with applications to wavefront propagation}. \emph{Stochastic Processes and their Applications}, {\bf 82}(1):23-52. \url{http://dx.doi.org/10.1016/S0304-4149(99)00003-4}

\bibitem{VeretLDPAvg2000}
Veretennikov A. 2000. \newblock{On large deviations for SDEs with small diffusion and averaging}. \emph{Stochastic Processes and their Applications}, {\bf 89}:69-79. \url{http://dx.doi.org/10.1016/S0304-4149(00)00013-2}

\bibitem{FWbook}
{Freidlin MI, Wentzell AD.} 2012. \emph{Random perturbations of dynamical systems}. Springer, 3rd ed.

\bibitem{Dykman}
{Dykman MI, Maloney CM, Smelyanskiy VN, Silverstein M.} 1998.
\newblock{Fluctuational phase-flip transitions in parametrically driven oscillators.} \emph{Phys. Rev. E}, 57(5):5202-5212. \url{http://dx.doi.org/10.1103/PhysRevE.57.5202}

\bibitem{JPolking}
{Polking JC.} 2003. \newblock{ODE software for matlab}. \url{http://math.rice.edu/~dfield}.


\bibitem{Day2}
{Day MV.} 1994. \newblock{Regularity of boundary quasi-potentials for planar systems.} \emph{Appl Math Optim}. {\bf 30}:79-101. \url{http://dx.doi.org/10.1007/BF01261992}

\bibitem{Day1}
{Day MV, Darden T.} 1985. \newblock{Some regularity results on the Ventcel-Freidlin quasi-potential function.} \emph{Appl Math Optim}. {\bf 13}:259-282. \url{http://dx.doi.org/10.1007/BF01442211}



\bibitem{YagasakiSIAM96} 
Yagasaki K. 1996. \newblock{The Melnikov theory for subharmonics and their bifurcations in forced oscillations}. \emph{SIAM J. Appl. Math.}, {\bf 56}(6):1720-1765. \url{http://www.jstor.org/stable/2102594}

\bibitem{namSow}
Sri Namachchivaya N, Sowers RB. 2002.
\newblock Rigorous stochastic averaging at a center with additive noise.
\newblock \emph{Meccanica}. {\bf 37}:85-114.
\url{http://dx.doi.org/10.1023/A:1019614613583}

\bibitem{KarTay}
Karlin S, Taylor HM. 1981. \emph{A second course in stochastic processes}. Academic Press.

\bibitem{EthierKurtz}
Ethier SN, Kurtz TG. 1986. \emph{Markov processes}. John Wiley \& Sons, New Jersey.

\bibitem{DayCharBoundary}
Day MV. 1990. Large deviations results for the exit problem with characteristic boundary. \emph{Journal of mathematical analysis and applications}, {\bf 147}:134-153.
 \url{http://dx.doi.org/10.1016/0022-247X(90)90389-W}

\bibitem{ByrdFriedman} Byrd PF, Friedman MD. 1971. \emph{Handbook of elliptic integrals for engineers and scientists}. Springer-Verlag. \url{http://dx.doi.org/10.1007/978-3-642-65138-0}

\bibitem{Olver1997} Olver FWJ. 1997. \emph{Asymptotics and special functions.} AK Peters Ltd. 

\bibitem{SpilioMultiscaleLDP} 
Spiliopoulos K. 2013. Large deviations and importance sampling for systems of slow-fast motion. \emph{Applied Mathematics {\&} Optimization}, {\bf 67}(1):123-161 \url{http://dx.doi.org/10.1007/s00245-012-9183-z}

\bibitem{SpilioHu} 
Hu W, Spiliopoulos K. 2015. Hypoelliptic multiscale Langevin diffusions: large deviations, invariant measures and small mass asymptotics. \emph{arxiv preprint}. \url{https://arxiv.org/abs/1506.06181}

\bibitem{FischerChiarini} Chiarini A, Fischer M. 2014. On large deviations for small noise Ito processes. \emph{Advances in Applied Probability} {\bf 46}(4):1126-1147. \url{https://projecteuclid.org/euclid.aap/1418396246}



\end{thebibliography}
\end{document}